\documentclass[headinclude=false,footinclude=false,paper=letter,11pt,oneside,draft=false,headings=normal,abstract=true,numbers=enddot]{scrartcl}
\areaset{5.6in}{8.5in}
\usepackage{pdfsync}

\usepackage[ps2pdf,hyperref,colorlinks=false,bookmarks=true,hyperfootnotes=false]{hyperref}

\usepackage{scrpage2}

\usepackage[latin2]{inputenc}
\usepackage[T1]{fontenc}

\usepackage{float}

\usepackage[leqno]{amsmath}
\usepackage{amsthm,latexsym,exscale}
\usepackage{mathtools}
\usepackage{faktor}

\usepackage{savesym}

\usepackage{dsfont}
\usepackage{marvosym}
\usepackage{manfnt}
\usepackage{ulem}
\savesymbol{Cross}
\usepackage{bbding}

\usepackage[nodayofweek]{datetime}

\usepackage{graphics,graphicx}
\usepackage[all,ps,cmtip]{xy}

\usepackage{pstricks,pst-node,pst-tree,pstricks-add,pst-grad}
\usepackage{pst-3dplot}

\usepackage{xspace}
\usepackage{enumerate}

\usepackage[footnotesize,format=hang,labelfont=sc,justification=justified,singlelinecheck=false,margin=3em]{caption}

\usepackage[multiple,hang]{footmisc}

\usepackage{tensor}

\usepackage{mathrsfs}

\usepackage{MnSymbol}

\usepackage[ps2pdf]{thumbpdf}

\newtheoremstyle{theorem}{\bigskipamount}{\bigskipamount}{\normalfont\normalsize\upshape}{0pt}%
	{\normalfont\normalsize\scshape}{:}{1ex}{}

\theoremstyle{theorem}

	\newtheorem{thm}{Theorem}

	\newtheorem*{thm*}{Theorem}

	\newtheorem{cor}[thm]{Corollary}
	\newtheorem{lem}[thm]{Lemma}
	\newtheorem{prop}[thm]{Proposition}

	\newtheorem*{cor*}{Corollary}
	\newtheorem*{lem*}{Lemma}
	\newtheorem*{prop*}{Proposition}

	\newtheorem{df}[thm]{Definition}
	\newtheorem{rem}[thm]{Remark}

	\newtheorem*{df*}{Definition}
	\newtheorem*{rem*}{Remark}

\newtheoremstyle{proof}{\bigskipamount}{\bigskipamount}{\normalfont\normalsize\upshape}{0pt}%
	{\normalfont\normalsize\scshape}{:}{1ex}{}

\theoremstyle{proof}

	\newtheorem*{pr}{Proof}

\newtheoremstyle{example}{}{}{\normalfont\normalsize\upshape}{0pt}%
	{\normalfont\normalsize\scshape}{:}{1ex}{}

\theoremstyle{example}

	\newtheorem*{ex*}{Example}
	
\newtheoremstyle{proofwo}{\bigskipamount}{\bigskipamount}{\normalfont\normalsize\upshape}{0pt}%
	{\normalfont\normalsize\scshape}{}{}{}

\theoremstyle{proofwo}

	\newtheorem*{pr*}{}
	
\newtheoremstyle{claim}{\medskipamount}{0ex}{\normalfont\normalsize\upshape}{0pt}%
	{\normalfont\normalsize\scshape}{:}{1ex}{}
	
\theoremstyle{claim}

	\newtheorem*{claim*}{Claim}
	
	\newtheorem*{obs*}{Observation}

\numberwithin{equation}{section}

\DeclareMathOperator{\rk}{rk\,}
\DeclareMathOperator{\minset}{Min}
\DeclareMathOperator{\comm}{Comm}
\DeclareMathOperator{\normalizer}{N}
\DeclareMathOperator{\centralizer}{C}

\DeclareMathOperator{\colim}{colim}
\DeclareMathOperator{\id}{id}
\DeclareMathOperator{\proj}{pr}
\DeclareMathOperator{\iso}{Iso}

\addtokomafont{disposition}{\rmfamily\scshape}

\mathtoolsset{showonlyrefs,showmanualtags}
\newtagform{brackets}{[}{]}

\renewcommand*{\eqref}[1]{(\refeq{#1})}

\begin{document}%
\bibliographystyle{amsalphacustom}%
%
%
\newcommand{\elaborate}[1][{}]{\marginpar{\begin{center} \Large \HandPencilLeft\;\;\; \end{center}\vspace*{-3ex}{\tiny #1}}}
\newcommand{\getitdone}[1][{}]{\marginpar{\begin{center}\Large \IroningI\;\;\; \end{center}\vspace*{-3ex}{\tiny #1}}}
\newcommand{\getitdoneI}[1][{}]{\marginpar{\begin{center}\Large \IroningI\;\;\; \end{center}\vspace*{-3ex}{\tiny #1}}}
\newcommand{\getitdoneII}[1][{}]{\marginpar{\begin{center}\Large \IroningII\;\;\; \end{center}\vspace*{-3ex}{\tiny #1}}}
\newcommand{\getitdoneIII}[1][{}]{\marginpar{\begin{center}\Large \IroningIII\;\;\; \end{center}\vspace*{-3ex}{\tiny #1}}}
\newcommand{\rewrite}[1][{}]{\marginpar{\begin{center}\Large \Handwash\;\;\; \end{center}\vspace*{-3ex}{\tiny #1}}}
\newcommand{\think}[1][{}]{\marginpar{\begin{center} \dbend\;\;\; \end{center}\vspace*{-0.5ex}{\tiny #1}}}
\newcommand{\noidea}[1][{}]{\marginpar{\begin{center}\LARGE \Biohazard\;\;\; \end{center}\vspace*{-0.5ex}{\tiny #1}}}
\newcommand{\spoof}{\marginpar{\mbox{\tiny Do I have to prove that?}}}
\newcommand{\dagegen}[1][{}]{\marginpar{\raisebox{-2ex}{\includegraphics[width=7ex,bb= 0 0 200 148]{pix/dagegen.eps}} \\ {\tiny #1}}}
\newcommand{\nicefrac}[2]{\ensuremath{\frac{#1}{#2}}}
\newcommand{\Alpha}{\ensuremath{\mathrm{A}}}
\newcommand{\Zeta}{\ensuremath{\mathrm{Z}}}
\newcommand{\Nu}{\ensuremath{\mathrm{N}}}
\newcommand{\Tau}{\ensuremath{\mathrm{T}}}
\newcommand{\NN}{\ensuremath{\mathds{N}}}
\newcommand{\BB}{\ensuremath{\mathds{B}}}
\newcommand{\BBo}{\ensuremath{\overset{\circ}{\mathds{B}}}}
\newcommand{\BBc}{\ensuremath{\overline{\mathds{B}}}}
\newcommand{\ZZ}{\ensuremath{\mathds{Z}}}
\newcommand{\QQ}{\ensuremath{\mathds{Q}}}
\newcommand{\RR}{\ensuremath{\mathds{R}}}
\newcommand{\CC}{\ensuremath{\mathds{C}}}
\newcommand{\PP}{\ensuremath{\mathds{P}}}
\newcommand{\HH}{\ensuremath{\mathds{H}}}
\newcommand{\KK}{\ensuremath{\mathds{K}}}
\newcommand{\LL}{\ensuremath{\mathds{L}}}
\newcommand{\II}{\ensuremath{\mathds{I}}}
\newcommand{\EE}{\ensuremath{\mathds{E}}}
\newcommand{\FF}{\ensuremath{\mathds{F}}}
\newcommand{\ident}{\ensuremath{\mathds{1}}}
\newcommand{\sphere}{\ensuremath{\mathds{S}}}
\newcommand{\TT}{\ensuremath{\mathds{T}}}
\newcommand{\cantor}{\ensuremath{\mathpzc{C}}}
\newcommand{\const}{\ensuremath{\text{\itshape const}}}
\newcommand{\GH}{G\hspace{-0.65ex}\raisebox{-0.69ex}{H}}
\newcommand{\eidi}{eigentlich diskontinuierlich}
\newcommand{\koko}{co-kompakt}
\newcommand{\vpz}{virtuell-polyzyklisch}
\newcommand{\PD}{\ensuremath{\mathrm{PD}}}
\newcommand{\group}[1]{\ensuremath{\normalfont{#1}}}
%
%
%
\newcommand{\gr}[1]{\ensuremath{\mathscr{#1}}}
\newcommand{\GE}{\gr{E}}
\newcommand{\GV}{\gr{V}}
\newcommand{\Q}{\gr{Q}}
\newcommand{\T}{\gr{T}}
\newcommand{\G}{\gr{G}}
\newcommand{\E}{\gr{E}}
\newcommand{\V}{\gr{V}}
\newcommand{\Vist}{\gr{V}_\mathrm{iST}}
\newcommand{\I}{\gr{I}}
%
%
\newcommand{\temp}[1]{\ensuremath{\mathcal{#1}}}
\newcommand{\Temp}{\temp{T}}
\newcommand{\W}{\temp{W}}
\newcommand{\Wo}{\overset{\circ}{\temp{W}}}
\renewcommand{\S}{\temp{S}}
\renewcommand{\L}{\temp{L}}
\newcommand{\Lp}{{\temp{L}\vphantom{'}}^+}
\newcommand{\Lip}{{\temp{L}'}^+}
\newcommand{\Liip}{{\temp{L}''}^+}
\newcommand{\Lkp}[1][k]{{\temp{L}\vphantom{'}}^{{\scriptscriptstyle{(#1)}}+}}
\newcommand{\Lm}{{\temp{L}\vphantom{'}}^-}
\newcommand{\Lim}{{\temp{L}'}^-}
\newcommand{\Liim}{{\temp{L}''}^-}
\newcommand{\Lk}[1][k]{{\temp{L}\vphantom{'}}^{\scriptscriptstyle{(#1)}}}
\newcommand{\Lkm}[1][k]{{\temp{L}\vphantom{'}}^{{\scriptscriptstyle{(#1)}}-}}
\newcommand{\Lpm}{{\temp{L}\vphantom{'}}^{\pm}}
\newcommand{\Lipm}{{\temp{L}'}^{\pm}}
\newcommand{\Liipm}{{\temp{L}''}^{\pm}}
\newcommand{\Lkpm}[1][k]{{\temp{L}\vphantom{'}}^{{\scriptscriptstyle{(#1)}}\pm}}
\newcommand{\prp}{{\lambda\vphantom{'}}^+}
\newcommand{\prm}{{\lambda\vphantom{'}}^-}
%
%
\newcommand{\shortedge}{\nicefrac{1}{\sqrt{2}}}
\newcommand{\longedge}{\ensuremath{\sqrt{2}}}
\newcommand{\sd}{\sim_{\varepsilon}} 
\newcommand{\cT}[2]{\ensuremath{(#1)_{T_{#2}}}}
\newcommand{\cS}[2]{\ensuremath{(#1)_{S_{#2}}}}
%
%
\newcommand{\dX}{d^{}_X}
\newcommand{\dT}{d^{}_{\Temp}}
\newcommand{\dtree}{d^{}_{\T}}
\newcommand{\dH}{d^{}_{\mathcal{H}}}
\newcommand{\vb}[1]{\ensuremath{\mathfrak{#1}}}
\newcommand{\virtbnd}{\partial^{}_{\infty}}
\newcommand{\tb}{\partial^{}_T}
\newcommand{\gd}[1]{\vec{#1}}
\newcommand{\cmsb}{\overset{.}{\sim}}
\newcommand{\St}[1]{\ensuremath{\overline{\operatorname{st}(#1)}}}
\newcommand{\act}{\ensuremath{\curvearrowright}}
\newcommand{\tca}{\ensuremath{\curvearrowleft}}
\newcommand{\rquotient}[2]{\ensuremath{\smash{\raisebox{0.3ex}{\ensuremath{#1}}/\raisebox{-0.3ex}{\ensuremath{#2}}}}}
\newcommand{\restr}[1]{\!\raisebox{-0.5ex}{\ensuremath{\bigm|}}_{\ensuremath{#1}}}
\newcommand{\cat}{\ensuremath{\mathrm{C\/A\/T}}}
\newcommand{\mls}{\ensuremath{\mathrm{M\/L\/S}}}
\renewcommand{\qed}{\hspace*{1cm}\hfill\rule{1.3ex}{1.3ex}}
\newcommand{\abs}[1]{\ensuremath{\left\lvert#1\right\rvert}}
\newcommand{\norm}[1]{\ensuremath{\left\lVert#1\right\rVert}}
\newcommand{\Lie}[2]{\ensuremath{\mathcal{L}_{#1}#2}}
\newcommand{\ball}[2]{\BB^{}_{#1}(#2)}
\newcommand{\cball}[2]{\overline{\BB}^{\vphantom{x}}_{#1}(#2)}
\newcommand{\nbhd}[2]{\ensuremath{\BB_{#1}(#2)}}
\newcommand{\scalar}[2]{\ensuremath{\left\langle #1 | #2 \right\rangle}}
\newcommand{\metr}[2]{\ensuremath{\left\langle #1 , #2 \right\rangle}}
\newcommand{\bra}[1]{\ensuremath{\left\langle #1 \right\rvert}}
\newcommand{\ket}[1]{\ensuremath{\left\lvert #1 \right\rangle}}
\newcommand{\Gl}{\ensuremath{\mathrm{Gl}}}
\newcommand{\End}[2]{\ensuremath{\mathrm{End}_{#1}(#2)}}
\newcommand{\del}{\ensuremath{\partial}}
\newcommand{\grad}{\smash{\raisebox{0,45ex}{\ensuremath{\bigtriangledown}}}}
%
%
\newcommand{\radiusscale}{N_1}
\newcommand{\wallnumscale}{C_1}
\newcommand{\wallnumshift}{C_2}
\newcommand{\odistmax}{\underline{C}}
\newcommand{\odistmin}{\overline{C}}
\newcommand{\anglemax}{\overline{\omega}}
\newcommand{\wallmax}{\radiusscale}%
\newcommand{\anglesummin}{(n+1)\anglemax}
%
%
%
\sloppy
%
\title{Virtual boundaries of Hadamard spaces with admissible actions of higher rank}

\author{Sebastian Grensing}

\date{\today}
\maketitle
\begin{abstract} \noindent 
Any discrete action of a group on a locally compact Hadamard space extends to a topological action on the virtual boundary.
\textsc{Croke} and \textsc{Kleiner} introduced a class of so-called admissible actions and associated geometric data which determine the topological conjugacy class of the boundary action. They also posed the question whether their results hold for a wider class of actions.

We show that, for the natural generalization, their question has to be answered in the negative: There is an admissible action of higher rank on a pair of Hadamard spaces with equivalent geometric data and an equivariant quasi-isometry which does not extend continuously to the virtual boundaries.
\end{abstract}
%
%
%
%
%
\section{Introduction}
\label{sec:intro}
We consider groups acting properly discontinuously and cocompactly by isometries on a Hadamard space, by which we refer to a complete, connected, and simply connected length space of non-positive curvature in the sense of \textsc{Alexandrov}.
As an isometry of a Hadamard space extends canonically to a homeomorphism of its virtual boundary, consisting of the asymptoty classes of geodesic rays endowed with the topology of uniform convergence on compact sets, any such action induces a topological action on the virtual boundary.

\textsc{Croke} and \textsc{Kleiner} introduced in \cite{crkl_geod} a class of so-called admissible actions on Hadamard spaces and associated invariants, the geometric data of the action, which determine the topological conjugacy class of the boundary action.
In particular, they show that any equivariant quasi-isometry extends canonically to a homeomorphism between the virtual boundaries, provided the geometric data of two admissible actions coincide up to scale.
From the graph of groups decomposition of the admissible group the authors derive an arrangement of convex subspaces in the Hadamard space, the edge and vertex spaces.
In a non-positively curved graph manifold these subspaces correspond to the universal covers of the Seifert pieces, resp.\ their boundary tori.
Employing this decomposition they show that any geodesic ray which is not asymptotic to some vertex space can be approximated, up to uniformly sublinear error, by the quasi-isometric image of a geodesic ray in a so-called \textit{template}, a model space for the sequence of edge spaces intersected by the geodesic ray.

This article addresses the question raised by \textsc{Croke} and \textsc{Kleiner} in \cite{crkl_geod} as to what extent the methods unravelled therein might be applicable to a wider class of actions.
Geometrically pertinent is the case of an admissible action of higher rank, i.e., virtually free abelian groups of higher rank being carried by the edges in the graph of groups decomposition, and the generalized notion of geometric data for such an action.
Despite the fact that most of the construction as in \cite{crkl_geod} is indeed attainable in our generalized case, quasi-isometries fail to relate geodesic rays in the Hadamard space to geodesic rays in a template as approximatively as in the low dimensional case.
Exploiting this repercussion, we show that, for the natural generalization, the aforementioned question has to be answered in the negative:
\begin{thm*}
There is an admissible action of higher rank on a pair of Hadamard spaces with equivalent generalized geometric data and an equivariant quasi-isometry which does not extend continuously to the virtual boundaries.
\end{thm*}
%
In the last section we briefly demonstrate how for admissible rank $3$ actions the equivalence of geometric data can be established in case an equivariant quasi-isometry extends continuously to the virtual boundaries.

\section{Preliminaries}
\label{sec:hd_sp}
The present section reviews some material on metric spaces of non-positive curvature, isometric group actions upon them and graph of groups.
Detailed expositions of this material can be found in \cite{ba_npc}, \cite{brhae_npc} or \cite{gr_metr_npc}, as well as \cite{se_trees}, while we rather have adopted the notation from \cite{didu_grgr} .

\subsection{Hadamard spaces}
\label{sec:qi}
A map $\Phi$ between metric spaces $(X,d^{}_X)$ and $(X',d^{}_{X'})$ is called an $(L,A)$-quasi-isometric embedding if for all $x,y \in X$
\begin{equation}
\frac{1}{L} d^{}_{X'} \bigl( \Phi(x) , \Phi(y) \bigr) - A \leqslant d^{}_X \bigl( x , y \bigr) \leqslant L \, d^{}_{X'} \bigl( \Phi(x) , \Phi(y) \bigr) + A 
\end{equation}
holds.
An $(L,A)$-quasi-isometric embedding is called an $(L,A)$-quasi-isometry if its image is $A$-dense in $X'$.
Every quasi-isometry $\Phi$ has a quasi-inverse; i.e., there is a quasi-isometry $\Psi$ such that $\Psi \circ \, \Phi$ lies within finite Hausdorff distance of the identity.
A quasi-geodesic is the quasi-isometric embedding of an interval.
The image of a quasi-geodesic does in general not lie within finite Hausdorff distance of a geodesic, yet each quasi-isometry induces a homeomorphism between the ends of the spaces.

\label{sec:npc}
By a Hadamard space we mean a complete, connected and simply connected length space of non-positive curvature in the sense of \textsc{Alexandrov}.
The metric of a Hadamard space is convex; therefore it is uniquely geodesic by the Cartan-Hadamard theorem.
A locally compact Hadamard space is proper by the Hopf-Rinow theorem.\nopagebreak

\label{sec:iso_grp_act}
\nopagebreak
Let $\Gamma \act X$ be a geometric action, i.e., properly discontinuous and cocompactly by isometries, on a locally compact Hadamard space $X$.
Then, by the \pagebreak \v{S}varc-Milnor lemma, $\Gamma$ is finitely generated and for any $x \in X$ the mapping $\Gamma \rightarrow X \; , \, \gamma \mapsto \gamma . x$ is an equivariant quasi-isometry with respect to a word metric on $\Gamma$.
In particular, all word metrics with respect to finite generating sets on $\Gamma$ lie in the same quasi-isometry class.
Furthermore, each isometry $\gamma \in \Gamma$ is semisimple; i.e., its displacement function $x \mapsto \delta_{\gamma} (x) = d_X ( \gamma . x , x )$ attains its minimum $\abs{\gamma\,}_X$.
The translation length $\abs{\,.\,}_X$ is invariant under conjugation ${}^{\alpha}\gamma = \alpha \gamma \alpha^{-1}$.
For each $\gamma \in \Gamma$ the minimal displacement set $\minset_X (\gamma) = \bigl\{ x \in X \bigm| \delta_{\gamma} (x) = \abs{\gamma\,}_X \bigr\}$ is nonempty and invariant under every isometry commuting with $\gamma$.
Since the displacement function $\delta_{\gamma}$ is continuous and convex, $\minset_X (\gamma)$ is closed and convex.
If $\gamma$ is hyperbolic, then $\minset_X (\gamma)$ is the union of the axes of $\gamma$ and splits as a metric product $Y \times \RR$; here $Y$ is a complete convex subspace of $X$ and the real fibres are the axes of $\gamma$.
Every isometry commuting with $\gamma$ is compatible with this product structure and splits as a product of an isometry of $Y$ and a translation of the real factor.
If $\Lambda$ is a free abelian subgroup of $\Gamma$, the following properties hold:
\begin{itemize}
\item $\Lambda$ is of finite rank, say $\rk \Lambda = k$.
\item $\smash{\minset^{}_X (\Lambda) = \bigcap\nolimits^{}_{\gamma \in \Lambda} \minset^{}_X (\gamma)}$ is nonempty and splits as a metric product $\smash{Y \times \EE^k}$ of a complete convex subspace $Y$ of $X$ and a Euclidean space of dimension $k$.
\item $\minset_X (\Lambda)$ is invariant under the action of the normalizer $\normalizer_{\Gamma} (\Lambda)$
\item The centralizer $\centralizer_{\Gamma} (\Lambda)$ has finite index in $\normalizer_{\Gamma} (\Lambda)$.
\item Every isometry $\alpha$ normalizing $\Lambda$ is compatible with the abovementioned splitting of the minimal displacement set; i.e., if $x = (\overline{x} , x_{\EE}) \in \minset_X (\Lambda)$ with respect to the product structure, there are isometries $\overline{\alpha} \in \iso (Y)$ and $\alpha_{\EE} \in \iso (\EE^k)$ such that $\alpha . (\overline{x} , x_{\EE}) = (\overline{\alpha} . \overline{x} , \alpha_{\EE} . x_{\EE})$.
\item $\alpha_{\EE}$ is a translation of $\EE^k$ for each $\alpha \in \centralizer_{\Gamma} (\Lambda)$.
\item $\Lambda$ acts trivially on $Y$ and cocompactly on the Euclidean factor; hence the quotient is a $k$-torus.
\item The induced action of $\rquotient{\normalizer_{\Gamma} (\Lambda)}{\Lambda}$ on $Y$ is properly discontinuous.
\end{itemize}

\begin{lem} \label{lem:N_coco}
Let $\Gamma$ be a group acting geometrically on a locally compact Hadamard space $X$ and $\Lambda$ a free abelian subgroup.
If $\zeta_1 , \ldots , \zeta_k$ denote free generators $\Lambda$, then, for each $R>0$, the centralizer $\centralizer_{\Gamma} (\Lambda)$ acts cocompactly on $\smash{M_R = \{ x \in X \mid \sum_{i \leqslant k} \delta_{\zeta_i} ( x ) \leqslant R \}}$.
\end{lem}

\begin{pr}
Since $M_R$ is proper, it suffices to show that for any sequence $x_n \in M_R$ there is a sequence $\xi_n \in \centralizer_{\Gamma} (\Lambda)$ such that a subsequence of $\xi_n . x_n$ converges.
By cocompactness of the $\Gamma$-action we may assume, without loss of generality, that there are $\gamma_n \in \Gamma$ such that $y_n = \gamma_n . x_n$ converges to an $y_{\infty}$ and $d_X ({}^{\gamma_n}\zeta_i . y_{\infty} , y_{\infty} ) \leqslant 1 + R$ for all $n \in \NN$ and $i \leqslant k$.
Since $\Gamma$ acts properly discontinuously, the set $\Sigma = \{ {}^{\gamma_n}\zeta_i \mid n \in \NN , i \leqslant k \}$ is finite.
Thus, possibly after replacing $n \mapsto ( {}^{\gamma_n}\zeta_1 , \ldots , {}^{\gamma_n}\zeta_k ) \in \Sigma^k$ by a subsequence, we can assume it to be constant.
Then $\xi^{}_n = \gamma^{-1}_1 \gamma^{}_n \in \centralizer^{}_{\Gamma} (\Lambda)$ for all $n \in \NN$, and $\xi_n . x_n$ converges to $\gamma^{-1}_1 . y^{}_{\infty} \in M^{}_R$.
\qed
\end{pr}
\pagebreak

\begin{cor}
$\centralizer_{\Gamma} (\Lambda)$ as well as $\normalizer_{\Gamma} (\Lambda)$ act cocompactly on $\minset_{X} (\Lambda)$and are, as a consequence, finitely generated.
\end{cor}

\begin{cor} \label{cor:min_mr_fin_hd}
If for some $R>0$ the set $M_R$ is nonempty, the Hausdorff distance between $M_R$ and $\minset_{X} (\Lambda)$ is finite.
\end{cor}

\label{sec:virt_bnd}
Two geodesic rays in a locally compact Hadamard space $X$ are called asymptotic if their Hausdorff distance is bounded.
The virtual boundary $\virtbnd X$ is the set of asymptoty classes of geodesic rays in $X$ endowed with the topology of uniform convergence on compact subsets.
For each $x \in X$ and $\vb{x} \in \virtbnd X$ there is a unique geodesic ray starting at $x$ and representing $\vb{x}$ which we denote by $\overline{x\vb{x}}$.
For a geodesic ray $c$ let $c(\infty)$ denote its equivalence class in $\virtbnd X$.
Every isometry of $X$ extends to a homeomorphism of its virtual boundary.
An isometric action $\Gamma \act X$ therefore induces a topological action $\Gamma \act \virtbnd X$.
Let $C$ be a closed convex subspace of $X$.
Then $\virtbnd C$ is closed in $\virtbnd X$.
If $C'$ is another such subspace within finite Hausdorff distance of $C$, then $\virtbnd C = \virtbnd C'$.
The virtual boundary of the product of two Hadamard spaces is naturally homeomorphic to the topological join of their virtual boundaries.
In particular, the virtual boundary of $X \times \EE^n$ is homeomorphic to the $n$-fold suspension of $\virtbnd X$.

\label{sec:hyp_sp}
Recall that an $(L,A)$-quasi-geodesic in a Hadamard space $X$ does not necessarily lie within finite Hausdorff distance of a geodesic.
But in the case that $X$ is $\delta$-hyperbolic, there even is a constant $C$, depending on $L,A$ and $\delta$ only, such that the image of any $(L,A)$-quasi-geodesic segment has Hausdorff distance at most $C$ to the geodesic segment connecting its endpoints.
Each quasi-geodesic ray thus defines a point in the virtual boundary, and every quasi-isometry between hyperbolic spaces extends continuously to a homeomorphism of their virtual boundaries.

\subsection{Graphs of groups}
\label{sec:gr_grp}
A graph of groups $\Gamma_{\G} (\_)$ consists of a connected graph $\G$ with edges $\E$ and vertices $\V$ together with a family of groups $\Gamma_{\G} (A)$, $A \in \E \cup \V$ and monomorphisms \linebreak $\iota^{\pm}_E \, : \, \Gamma_{\G} (E) \hookrightarrow \Gamma_{\G} (\del^{\pm} E)$ for each $E \in \E$.
An action of a group $\Gamma$ on a connected graph $\gr{X}$ induces a graph of groups on the quotient $\G = \rquotient{\gr{X}}{\Gamma}$ as follows:
For each edge or vertex $A$ of $\G$ fix an $a \in A$ and let $\Gamma^{}_{\G} (A) = \Gamma^{}_a$.
For the representative $e$ of an edge $E$ let $\gamma^{\pm}_e \in \Gamma$ such that $\gamma^{\pm}_e \del^{\pm} e$ represents $\del^{\pm} E$. 
For an edge $E$ of $\G$ and the representatives $e$ and $v^{\pm}$ of $E$ and $\del^{\pm} E$ the required monomorphisms are given by $\alpha \mapsto {}^{\gamma^{\pm}_e} \alpha = \gamma^{\pm}_e \alpha {\gamma^{\pm}_e}^{-1}$ from $\Gamma_{\G} (E) = \Gamma_e < \Gamma_{\del^{\pm} e}$ to ${}^{\gamma^{\pm}_e} \Gamma_{\del^{\pm} e} = \Gamma_{\gamma^{\pm}_e . \del^{\pm} e} = \Gamma_{\G} (\del^{\pm} E)$.

\label{sec:fund_grp}
Let $\Gamma_{\G} (\_)$ be a graph of groups and $\G_0$ a maximal subtree of $\G$ with edges $\E_0$.
The fundamental group%
\footnote{Let $\G'$ denote the barycentric subdivision of $\G$. A graph of groups $\Gamma_{\G} (\_)$ is then a $\G'$-diagram $\mathsf{G}$ in the category of groups with monomorphisms, and $\pi_1 (\Gamma_{\G} (\_) , \G_0) \cong \colim \mathsf{G}$.}
$\pi_1 (\Gamma_{\G} (\_) , \G_0)$ of $\Gamma_{\G} (\_)$ at $\G_0$ is defined by the presentation given by the generators $\bigl\{ \tau_E \mid E \in \E \bigr\} \cup \bigcupdot\nolimits_{V \in \V} \Gamma_{\G} (V)$ and the relations of each vertex group together with $\tau^{}_E \iota^-_E (\gamma) \tau^{-1}_E = \iota^+_E (\gamma)$ for all edges $E$ of $\G$ and $\gamma \in \Gamma_{\G} (E)$ as well as $\tau^{}_E = 1$ if $E \in \E_0$.
The isomorphism class of the fundamental group depends on neither the choice of the maximal subtree $\G^{}_0$ nor the isomorphism class of the graph of groups.

\begin{ex*}
Let $\G$ be the graph with a single edge $E$ incident to distinct vertices $V$, $W$ \linebreak and $\Gamma^{}_{\G} (\_)$ a graph of groups with $A = \Gamma^{}_{\G} (V)$, $B = \Gamma^{}_{\G} (W)$ and $C = \Gamma^{}_{\G} (E)$.
Its fundamental group is isomorphic to the amalgamated free product%
\footnote{In this case, $\G' = \bullet \leftarrow \bullet \rightarrow \bullet$ and a graph of groups over $\G$ is a pushout diagram of two monomorphisms.}
$A \ast^{}_C B$.
\end{ex*}

\begin{ex*}
Let $\G$ be the graph with one edge $E$ and one vertex $V = \del^- E = \del^+ E$ only.
A graph of groups over $\G$ then consists of two groups $A$ and $C$ and two monomorphisms from $C$ into $A$.
Its fundamental group%
\footnote{Cutting a Riemannian surface $X$ of positive genus along a non-separating simple closed geodesic results in a Riemannian surface $X'$ with two additional boundary components $Y^{}_1$ and $Y^{}_2$. The embeddings $\sphere^1 \simeq Y^{}_i \hookrightarrow X'$ induce two monomorphisms from $\ZZ$ into $\pi^{}_1 (X')$ and $\pi^{}_1 (X) = \pi^{}_1 (X') \ast^{}_{\ZZ}$.}
is an HNN-extension $A \ast^{}_C$ of $A$ over $C$.
\end{ex*}

The fundamental group of a graph of groups $\Gamma^{}_{\G} (\_)$ over a finite graph $\G$ at a maximal subtree $\G^{}_0 \subset \G$ can be calculated by employing the two basic constructions mentioned above:
For each edge $E$ of $\G^{}_0$ a new graph of groups is defined by replacing $E$ together with its adjacent vertices by a single vertex carrying the amalgamated free product $\smash{\Gamma^{}_{\G} (\del^- E) \ast_{\Gamma^{}_{\G} (E)} \Gamma^{}_{\G} (\del^+ E)}$.
Finitely many such reductions yield a rose shaped graph, each petal corresponding to an edge of $\G \setminus \G^{}_0$.
The fundamental group now results from successively deleting the remaining edges while replacing the single vertex group by the respective HNN-extension.

On the other hand, for a graph of groups $\Gamma_{\G} (\_)$ with fundamental group $\Gamma$ the sets $\widetilde{\V} = \bigcupdot\nolimits_{V \in \V} \rquotient{\Gamma}{\Gamma^{}_{\G} (V)}$ and $\widetilde{\E} = \bigcupdot\nolimits_{E \in \E} \rquotient{\Gamma}{\Gamma^{}_{\G} (E)}$ together with the boundary maps $\gamma \Gamma^{}_{\G} (E)  \longmapsto \gamma \Gamma^{}_{\G} (\del^{\pm} E)$ define a tree $\T$ upon which $\Gamma$ acts; it is called its Bass-Serre tree.
The graph of groups induced on the quotient $\rquotient{\T\!}{\Gamma}$ is isomorphic to $\Gamma_{\G} (\_)$.
For a maximal subtree $\G^{}_0 \subset \G$ each generator $\tau_E$ of an edge $E$ not in $\G^{}_0$ corresponds to an element of the topological fundamental group of $\G$ and acts as a translation on $\T$; in particular, $\tau_E . \del^- e = \del^+ e$ holds for any lift $e$ of $E$.

\section{Admissible actions}
In what follows, we will generalize the class of group actions which \textsc{Croke} and \textsc{Kleiner} introduced in \cite{crkl_geod} as admissible actions.
Therein the authors classified admissible actions by a family of functions defined on the edge and vertex groups.
Additionally there is an evident generalization of these so-called geometric data for the wider class of actions we consider here.

\subsection{Admissible graphs of groups}
Recall that two subgroups $A$ and $B$ of a group $G$ are commensurable if $A \cap B$ has finite index in both $A$ and $B$; we will then write $A \sim B$.
They are conjugate commensurable if there is a $g \in G$ such that ${}^g A = \{ g a g^{-1} \mid a \in A \} \sim B$.
Let $\comm_G (A) = \{ g \in G \mid {}^g A \sim A \}$ denote the commensurizer subgroup of $A$ in $G$.

\begin{df} \label{def:adm_gr_gr}
A graph of groups $\Gamma^{}_{\G} (\_)$ is admissible of rank $k$ if $\G$ is finite, contains at least one edge and the following conditions are satisfied:
\begin{enumerate}
	\item Each vertex group $\Gamma^{}_{\G} (V)$ is an extension of a non-elementary hyperbolic group by a free abelian subgroup $\Lambda (V) < \Gamma^{}_{\G} (V)$ of rank $k-1$.
	\label{def:adm_gr_gr:vertex_groups}
	
	\item Each edge group $\Gamma^{}_{\G} (E)$ is vitually free abelian of rank $k$.
	\label{def:adm_gr_gr:edge_groups}
	
	\item If $E$ is an edge incident to $V = \del^{\epsilon} E$, then for every $\gamma \in \Gamma^{}_{\G} (V) \setminus \iota^{\epsilon}_E \left( \Gamma^{}_{\G} (E) \right)$ the image of $\Gamma^{}_{\G} (E)$ in $\Gamma^{}_{\G} (V)$ and its conjugate under $\gamma$ are incommensurable. For another edge $E'$ incident to $V$ the images of $\Gamma^{}_{\G} (E)$ and $\Gamma^{}_{\G} (E')$ in $\Gamma^{}_{\G} (V)$ are not conjugate commensurable.
	\label{def:adm_gr_gr:non_comm}
	
	\item For any edge $E$ the preimages%
	\footnote{By \ref{def:adm_gr_gr:vertex_groups} we can assume each subgroup $\Lambda (\partial^{\pm} E)$ to lie in the image of ${\iota^{\pm}_E}$ since a finite extension of a hyperbolic group is hyperbolic as well.}
	${\iota^-_E}^{-1} \Lambda (\partial^- E)$ and ${\iota^+_E}^{-1} \Lambda (\partial^+ E)$ generate a subgroup of finite index in $\Gamma_{\G} (E)$.
	\label{def:adm_gr_gr:vertex_union}
\end{enumerate}

A geometric action of a group $\Gamma$ on a locally compact Hadamard space is admissible of rank $k$ if $\Gamma$ is isomorphic to the fundamental group of an admissible graph of groups of rank $k$.
\end{df}

\begin{ex*}[Torus complexes] \label{ex:adm_tc}
Let $T_i$, $i \leqslant m+1$, be a family of flat $2$-Tori, each containing a pair of non-homotopic simple closed geodesics $\alpha_i$ and $\beta_i$ such that \linebreak $l(\beta_i) = l(\alpha_{i+1})$.
Identifying $\beta_i$ and $\alpha_{i+1}$ and endowing the resulting space with the induces length metric yields a compact space $Y$ of non-positive curvature.

This construction determines a graph with $m$ vertices $V_i = \rquotient{( T_i \cup T_{i+1})}{\{ \beta_i = \alpha_{i+1} \!\} }$ and $m-1$ edges $E_i = V_i \cap V_{i+1}$.
The resulting graph of groups, each edge or vertex carrying its fundamental group, is admissible of rank $2$:
Each vertex is homeomorphic to the product of a figure eight and a loop; hence its fundamental group is isomorphic to $\FF_2 \times \ZZ$.
Every edge $E_i$ is incompressible in each of the vertices incident to it, the maps from $\pi_1 (E_i)$ in $\pi_1 (V_i)$ and $\pi_1 (V_{i+1})$ are therefore injective.
The generators $\alpha_i, \beta_i$ of the edge group $\pi_1 (E_i) \simeq \ZZ^2$ each map to a generator of either the infinite cyclic or the free factor of an adjacent vertex group.

The universal cover $X$ of $Y$ is, by the Cartan-Hadamard theorem, a locally compact Hadamard space upon which $\pi_1 (Y)$ acts geometrically.
By the van Kampen theorem $\pi_1 (Y)$ is isomorphic to the fundamental group of the abovementioned graph of groups.
The action $\pi_1 (Y) \act X$ is therefore admissible of rank $2$.
\end{ex*}

\begin{ex*}[Generalized graph manifolds] \label{ex:adm_gr_mf}
A generalized graph manifold of dimension $n+2$ is a closed differentiable manifold $M$ consisting of a finite union of blocks \linebreak $M_i \cong N_i \times \TT^n$, where each $N_i$ is a compact, orientable Riemannian surface $N_i$ with nontrivial boundary and of negative Euler characteristic.
The connected components of the boundary of a block are $(n+1)$-tori, and blocks are glued amongst themselves by identifying two boundary tori via diffeomorphisms in such a way that no two torus fibres of different blocks are homotopic.

This defines a graph with the blocks as vertices and an edge between two vertices if connected components of their boundary are identified.
Letting each edge or vertex carry the fundamental group of its defining block or boundary component yields a graph of groups whose edge and vertex groups are free abelian of rank $n+1$ and direct products of a fuchsian group and a free abelian group of rank $n$ respectively.

In case $M$ carries a metric of non-positive sectional curvature%
\footnote{See \cite{svet_obstr} for necessary and sufficient conditions for the existence of such a metric.},
it follows, just like in the previous example, that the action of the fundamental group on the universal cover on $M$ is admissible of rank $n+1$.
\end{ex*}

\begin{ex*}[Poincar\'e duality groups] \label{ex:adm_pd_grp}
If a group is isomorphic to the fundamental group of a graph of groups, then each edge represents a splitting of the group as an amalgamated product or an HNN-extension.
A graph of group decomposition of a group which thereby comprises all possible splittings as amalgamated product or HNN-extension over a certain class of groups up to conjugacy is called JSJ-decomposition%
\footnote{This term has been coined by \textsc{Sela} in \cite{se_jsj} for the characterization of graph of groups decompositions of torsion-free hyperbolic groups, referring to the work by \textsc{Waldhausen} \cite{wald_3mf} on characteristic submanifolds of $3$-manifolds and the subsequent results by \textsc{Jaco}, \textsc{Shalen} \cite{jash_seifert} and \textsc{Johannson} \cite{joh_3mf}.}%
.
JSJ-decompositions of Poincar\'e duality groups over virtually polycyclic groups are presented in \cite{scsw_reg_nbhd} and \cite{scsw_at}.

Let $\Gamma$ be a Poincar\'e duality group of dimension $n+1$ acting geometrically on a locally compact Hadamard space $X$ and assume the edges of special Seifert type and the atoroidal edges of its JSJ-decomposition as in \cite[Theorem 5.1]{scsw_at} to be trivial.
Then the action of $\Gamma$ on $X$ is admissible of rank $n$.

\textsc{Wall} conjectured in \cite[Conjecture 10.5, p. 91]{wall_geom_groups} that strictly atoroidal Poincar\'e duality pairs of dimension $3$ are in fact relative hyperbolic.
Should this prove true, one should be able to adopt techniques from \cite{hu_iso_flats} and \cite{hukl_iso_flats} in order to obtain results on actions of Poincar\'e duality groups of dimension $3$ on Hadamard spaces without the aforementioned restriction.
\end{ex*}

\begin{lem} \label{rem:adm_tr_gr}
Let $\Gamma$ be a group which acts geometrically on a locally compact Hadamard space $X$.
This action is admissible of rank $k$ if and only if there is a $\Gamma$-finite $\Gamma$-tree $\T$ and a $\Gamma$-equivariant family $\Lambda^{}_a$ of normal subgroups of the stabilizers $\Gamma^{}_a$ of the simplices $a$ of $\T$ with the following properties:
\begin{enumerate}
\item For each vertex $v$ the subgroup $\Lambda^{}_v$ is free abelian of rank $k-1$, $\comm^{}_{\Gamma} (\Lambda^{}_v) = \Gamma^{}_v$ and $\rquotient{\Gamma^{}_v}{\Lambda^{}_v}$ is non-elementary hyperbolic.
	\label{rem:adm_tr_gr:vertex_groups}
	
\item For each edge $e$ the subgroup $\Lambda^{}_e$ is free abelian of rank $k$, $\rquotient{\Gamma^{}_e}{\Lambda^{}_e}$ is finite and $\Lambda^{}_{\partial^{\pm} e} < \Lambda^{}_e$.
	\label{rem:adm_tr_gr:edge_groups}

\item For any two consecutive edges $e$ and $e'$ the subgroups $\Lambda^{}_e$ and $\Lambda^{}_{e'}$ are incommensurable.
	\label{rem:adm_tr_gr:non_comm}
\end{enumerate}
\end{lem}

\begin{pr}
To begin with, we will show that the Bass-Serre tree $\T$ of an admissible graph of groups satisfies the conditions above:
For each edge or vertex $a$ of $\T$ let $\Lambda^{}_a$ be the equivariant lift of the normal free abelian subgroup of rank $k$ or $k-1$ resp.\ in $\Gamma_{\G} (\Gamma . a)$.
Let $e$ be an edge of $\T$ and $v$ a vertex incident to it.
Since $\Gamma^{}_v$ acts by semi simple isometries on the Hadamard space $X$ and $\Lambda^{}_v$ is normal in $\Gamma^{}_v$, a finite index subgroup of $\Gamma^{}_v$ centralizes $\Lambda^{}_v$.
Thus, for any $\xi \in \Lambda^{}_v$
\begin{equation}
\Gamma^{}_{\xi . e} = {}^{\xi} \Gamma^{}_e \sim {}^{\xi} \bigl( \Gamma^{}_e \cap \centralizer^{}_{\Gamma^{}_v} (\Lambda^{}_v) \bigr) = \Gamma^{}_e \cap \centralizer^{}_{\Gamma^{}_v} (\Lambda^{}_v) \sim \Gamma^{}_e
\end{equation}
holds, and $\xi \in \Gamma^{}_e$ by Definition \ref{def:adm_gr_gr}.
Since $\G$ is finite, we can assume, possibly after replacing $\Lambda^{}_v$ by a characteristic subgroup of finite index, that $\Lambda^{}_v < \Lambda^{}_e$ for any vertex $v$ incident to an edge $e$.
It remains to show that $\comm^{}_{\Gamma} (\Lambda^{}_v) < \Gamma^{}_v$ for every vertex $v$.
Note that for any vertices $v$ and $w$ adjacent to the same edge both $\Lambda^{}_v \cap \Lambda^{}_w$ and $\Lambda^{}_v \Lambda^{}_w$ are free abelian.
$\Lambda^{}_v \Lambda^{}_w$ has finite index in $\Lambda^{}_e$ by Definition \ref{def:adm_gr_gr} \ref{def:adm_gr_gr:vertex_union} and is therefore of rank $k$ by Definition \ref{def:adm_gr_gr} \ref{def:adm_gr_gr:edge_groups}; hence $\Lambda^{}_v \cap \Lambda^{}_w$ is of rank $k-2$.
For any edges $e \neq e'$ incident to the same vertex $v$ the subgroups $\Gamma^{}_e \cap \Gamma^{}_{e'}$ and $\Lambda^{}_v$ are commensurable.
Fix $\gamma \in \Gamma \setminus \Gamma^{}_v$.
There is a geodesic path $(v^{}_0 e^{}_1 v^{}_1 \ldots e^{}_n v^{}_n)$ of length $n \geqslant 3$ in $\T$ with $v^{}_1 = v$ and $v^{}_{n-1} = \gamma . v$.
Since $\Lambda^{}_v < \Lambda^{}_{e^{}_1} < \Gamma^{}_{v^{}_0}$ and ${}^{\gamma} \Lambda^{}_v = \Lambda^{}_{\gamma . v} < \Gamma^{}_{v^{}_n}$, we deduce from
\begin{equation}
\Lambda^{}_v \cap {}^{\gamma} \Lambda^{}_v < \Gamma^{}_{v^{}_0} \cap \Gamma^{}_{v^{}_n} < ( \Gamma^{}_{e^{}_1} \cap \Gamma^{}_{e^{}_2} ) \cap ( \Gamma^{}_{e^{}_2} \cap \Gamma^{}_{e^{}_3} ) \sim \Lambda^{}_{v^{}_1} \cap \Lambda^{}_{v^{}_2} \cong \ZZ^{k-2}
\end{equation}
that $\Lambda^{}_v$ and ${}^{\gamma} \Lambda^{}_v$ are not commensurable.

Now, let $\T$ be a tree and $\Gamma \act \T$ an action with the required properties.
The graph of groups on the quotient $\G = \rquotient{\T}{\Gamma}$ satisfies the first three conditions of Definition \ref{def:adm_gr_gr} quite obviously.
Let $E$ be an edge of $\G$ and $e \in E$ a lift to $\E\!\T$.
If $\Lambda^{}_{\partial^- e} \Lambda^{}_{\partial^+ e}$ were not of finite index in $\Lambda^{}_e$, the subgroups $\Lambda^{}_{\partial^- e}$ and $\Lambda^{}_{\partial^+ e}$ were commensurable and $\Gamma^{}_{\partial^+ e} < \comm_{\Gamma} (\Lambda^{}_{\partial^- e})$ contradicted property \ref{rem:adm_tr_gr:vertex_groups}.
\qed
\end{pr}

In particular we have for an admissible action $\Gamma \act X$ of rank $k$ that $\centralizer_{\Gamma} (\Lambda_a) \sim \normalizer_{\Gamma} (\Lambda_a) = \comm_{\Gamma} (\Lambda_a) = \Gamma_a$ for any simplex $a$ of the Bass-Serre tree of the graph of groups decomposition of $\Gamma$.

\begin{rem} \label{rem:adm_tr_gr:uniqueness}
Let $\Gamma$ be a group, $X$ a locally compact Hadamard space and $\Gamma \act X$ an admissible action of rank $k$. 
The decomposition of $\Gamma$ as an admissible graph of groups is unique up to isomorphism.
\end{rem}

We omit the proof, which is an easy variation on \cite[Lemma 3.7]{crkl_geod}.

\subsection{Geometric data}
\label{sec:geom_data_def}
We will now define the geometric data of an admissible action.
Let $\Gamma \act X$ be admissible of rank $n+1$ and let $\T$ denote the Bass-Serre tree of the graph of groups decomposition of $\Gamma$.

Recall, that for a vertex $v$ of $\T$ the minimal displacement set $Y_v = \minset_X (\Lambda_v)$ is a nonempty, $\Gamma_v$-invariant, closed and convex subspace of $X$.
It splits as a metric product $Y_v \cong \overline{Y}_v \times \EE_v$ of a convex subspace $\overline{Y}_v$ and an $n$-dimensional Euclidean space $\EE_v$.
The isometries in $\Gamma_v$ respect the product structure in the sense that for any $\gamma \in \Gamma_v$ there are isometries $\overline{\gamma} \in \iso (\overline{Y}_v)$ and $\gamma_{\EE} \in \iso (\EE_v)$ such that $\gamma . (\overline{x} , x_{\EE}) = (\overline{\gamma} . \overline{x} , \gamma_{\EE} . x_{\EE})$ for any $x \in Y_v$ with coordinates $(\overline{x} , x_{\EE})$ with respect to the splitting.
Furthermore, the action of $\Lambda_v$ is trivial on the first factor and as a lattice on the Euclidean factor.
The action of $\Gamma_v$ on $\overline{Y}_v$ therefore factorizes over $\Lambda_v$, and $\rquotient{\centralizer_{\Gamma} (\Lambda_v)}{\Lambda_v}$ acts geometrically on the locally compact Hadamard space $\overline{Y}^{}_v$.
In particular, $\overline{Y}^{}_v$ is a hyperbolic space with more than two ends by Lemma \ref{rem:adm_tr_gr} \ref{rem:adm_tr_gr:vertex_groups}.

The first geometric datum at a vertex $v$ is defined, analogous to \cite[Definition 3.9]{crkl_geod}, as the function on $\Gamma_v$ which maps each isometry to its minimal translation length on the hyperbolic factor of $Y_v$:
\begin{equation} \label{def:mls}
\mls_v \, : \, \Gamma_v \longrightarrow \RR_{\geqslant 0} \; , \, \gamma \longmapsto \abs{\overline{\gamma}\,}_{\overline{Y}_v}
\end{equation}
Since the action of $\Gamma_v$ on $\overline{Y}_v$ factors over $\Lambda_v$, we may also regard to $\mls_v$ as being defined on their hyperbolic quotient%
\footnote{In case $X$ is a generalized graph manifold, the map $\mls_v$ corresponds to the marked length spectrum of the Riemannian surface in the block $v$.}.

In order to define the second geometric datum we consider for a vertex $v$ the real vector space $\Lambda_v \otimes_{\ZZ} \RR \cong \RR^n$.
The mapping $\Lambda_v \rightarrow \RR \; , \zeta \mapsto \smash{\abs{\zeta}^2_X}$ induces a scalar product $<.\,,.>$ on $\Lambda_v \otimes_{\ZZ} \RR$ whose norm we denote by $\smash{\norm{\,.\,}}$.
Thus, by virtue of the $\Lambda_v$-action, the Euclidean factor $\EE_v$ of $Y_v$ becomes canonically isometric to $\Lambda_v \otimes^{}_{\ZZ} \RR$ and a choice of generators $\zeta^{}_1 , \ldots , \zeta_n$ of $\Lambda_v$ yields a basis of $\Lambda_v \otimes_{\ZZ} \RR$.
The centralizer of $\Lambda_v$ acts by translations on $\EE_v$.
For each $\gamma \in \centralizer_{\Gamma} (\Lambda_v)$ let $\vec{\gamma} = ( g_1 , \ldots , g_n )$ denote the coordinates of its translation vector $\gamma_{\EE}$ in $\Lambda_v \otimes_{\ZZ} \RR \cong \EE_v$ with respect to the basis $\zeta_1 , \ldots , \zeta_n$.
We now choose, in an equivariant way, generators of $\Lambda_v$ for each vertex $v$ in $\T$.
The homomorphism of groups%
\begin{equation} \label{def:tau}
\Tau^{}_v \, : \, \centralizer^{}_{\Gamma} (\Lambda^{}_v) \longrightarrow \RR^n \; , \, \gamma \longmapsto \vec{\gamma}
\end{equation}
then constitutes the second geometric datum at the vertex $v$.
The following example illustrates how this mapping does indeed carry geometric information, even if $\Gamma_v$ is a product.

\begin{ex*}
Let $Y$ be the Cayley graph of the free group $\FF_2$ on two generators $f_1 , f_2$ and $\phi$ the homomorphism from $\FF_2$ to $\ZZ$ defined by $f_1 \mapsto 1$ and $f_2 \mapsto 0$.
For each $t \in \RR$ the group $\FF_2 \times \ZZ$ acts on $Y \times \EE$ by $(g,k) . (y,r) = (g . y , r + k + t \phi(g) )$.
Choosing $1$ as generator of $\ZZ$ we then have $\Tau (g,k) = k + t \phi(g)$ for $(g,k) \in \FF_2 \times \ZZ$.
\end{ex*}

The families $\{ \mls_v \}_{v \in \V}$ and $\{ \Tau_v \}_{v \in \V}$ are called the (generalized) geometric data of the admissible action $\Gamma \act X$.
They are invariant under the action of $\Gamma$, for the families $\{ Y_v \}_{v \in \V}$ are equivariant and the translation length of an isometry is invariant under conjugation.
In particular, the geometric data can be considered as being defined on the vertex groups of the graph of groups over $\G = \rquotient{\T}{\Gamma}$.

Let $\Gamma \act X'$ be another admissible action of rank $n+1$ and $\{ \mls'_v \}^{}_{v \in \V}$ together with $\{ \Tau'_v \}^{}_{v \in \V}$ its geometric data.
The geometric data of the two actions are equivalent if for each vertex $v$ there is a $\lambda (v) \in \RR^{}_{>0}$ such that $\mls'_v = \lambda (v) \mls^{}_v$ and $\Tau'_v = \Tau^{}_v$ hold%
\footnote{Presumably, the choices of generators of $\Lambda_v$ coincide for both actions.}.

\section{Proof of the Theorem}
\label{sec:cex}
In this section we will construct a pair of admissible rank $3$ actions with equivalent geometric data and an equivariant quasi-isometry between the respective Hadamard spaces which does not extend continuously to a homeomorphism of the virtual boundaries.

The perception of this example arose from our efforts to obtain a generalization of Theorem 1.3 in \cite{crkl_geod} for admissible actions of higher rank.
In their paper \textsc{Croke} and \textsc{Kleiner} define the edge and vertex spaces, mentioned in the introduction, to be certain tubular neighbourhoods of the minimal displacement sets of the free abelian subgroups of the edge and vertex stabilizers resp.; even though their construction carries over nicely to our generalized case we will refrain from making this precise.
The authors' proof relies heavily on the fact that in case of a rank 2 admissible action any geodesic segment in the Hadamard space diverges at most sublinearly from the quasi-isometric image of a corresponding geodesic segment in the so-called \textit{template} of its route; a template is a piecewise Euclidean model for the sequence of edge and vertex spaces a geodesic segment or ray traverses.
Given an admissible action of higher rank there are geodesic rays that pass through a sequence of edge and vertex spaces with a common Euclidean factor of positive dimension.
The corresponding template therefore decomposes into a product of an Euclidean factor and a template of lower dimension.
Even though the unit speed reparametrization of the projection of any segment of such a ray onto the non-Euclidean factor corresponds up to sublinear error to a geodesic segment in the lower dimensional template, still, the geodesic segment and the quasi-isometric image from the template might diverge at a rate proportional to arclength.
By gluing such products we are able to construct a sequence of geodesic segments that violate any sublinear estimate as in \cite[Theorem 5.1]{crkl_geod}.

\subsection{Construction}
\label{sec:cex:construction}
First of all, we define a continuous family of homotopically equivalent compact, non-positively curved length spaces $Y^{}_{\varepsilon}$.
The actions of $\Gamma = \pi^{}_1 (Y^{}_{\varepsilon})$ on their universal covers will then be shown to be admissible of rank 3 with equivalent geometric data.

The topology of $Y^{}_{\varepsilon}$ is easy to describe: We fix a disjoint family of nine $3$-tori with a choice of pairs of non-homotopic incompressible $2$-tori in each of them as well as homeomorphisms between pairs of these $2$-tori lying in distinct $3$-tori. 
Glueing the $3$-tori along the mapping cylinders of these homeomorphisms yields a connected space.
The thorough construction goes as follows:

For $\varepsilon \geqslant 0$ let $\varepsilon^{}_i = \varepsilon$ for $i = 2,3$ and $\varepsilon^{}_i = \varepsilon^{}_{\times} = 0$ for $i = 0,1$.
Now define $E^{}_{\times} \cong E^{\pm}_i \cong \TT^3 \times [0, \varepsilon_i]$ to be nine (thickened) flat tori with coordinates $(r,s,t;u)$ and radii given \nolinebreak by
\begin{align}
	\qquad (r,s,t;u)^{\pm}_i & = (r + \shortedge k,s + \shortedge l,t + \longedge m;u)^{\pm}_i
		& \text{ for } i = 2,3 \; , \qquad  \\
	(r,s,t;0)^{\pm}_i & = (r + k,s + l,t + \longedge m ; 0)^{\pm}_i
		& \text{ for } i = 0 , 1 \; ,  \qquad \\
\shortintertext{and}
	(r,s,t;0)^{}_{\times} & = (r + \longedge k,s + \longedge l,t + \longedge m ; 0)^{}_{\times}
\end{align}
for any $( k , l , m ) \in \ZZ^3$.
These pieces are glued along ten cylinders $C^{\pm}_{\times} \cong C^{\pm}_i \cong \TT^2 \times [0,1]$ with coordinates $(r,s,u)^{\pm}_i$ and radii as to make the equivalence relation generated by the following equations become compatible with the $\ZZ^3$-action (compare Figure \ref{fig:Y0-_domains} for an illustration of the fundamental domains and the respective identifications for a two-dimensional subspace in case $\varepsilon = 0$):
\begin{align}
	(s,-s,t;0)^{\pm}_0 & \sim (s,t;0)^{\pm}_0
		& (s,t;1)^{\pm}_0 &\sim (s,-s,t;0)^{\pm}_1 \\
	(s,0,t;0)^{\pm}_1 & \sim (s,t;0)^{\pm}_1
		& (s,t;1)^{\pm}_1 & \sim (s,s,t;0)^{\pm}_2 \\
	(0,s,t;\varepsilon)^{\pm}_2 & \sim (s,t;0)^{\pm}_2
		& (s,t;1)^{\pm}_2 & \sim (0,s,t;\varepsilon)^{\pm}_3 \\
	(s,s,t;0)^{\pm}_3 & \sim (s,t;0)^{\pm}_3
		& (s,t;1)^{\pm}_3 & \sim (s,0,t;0)^{\pm}_0 \label{eq:Yd_gluing_koord} \\ \\
	(s,-s,t;0)^-_0 & \sim (s,t;0)^-_{\times}
		& (s,t;1)^-_{\times} & \sim (s,0,t;0)^{}_{\times} \\
	(s,t,0;0)^{}_{\times} & \sim (s,t;0)^+_{\times}
		& (s,t;1)^+_{\times} & \sim (-s,s,t;0)^+_0
\end{align}
Let $\alpha^{\pm}_i$, $\beta^{\pm}_i$, $\gamma^{\pm}_i$, and $\alpha^{}_{\times}$, $\gamma^+_{\times}$, $\gamma^-_{\times}$ denote the loops defined by the coordinate lines on the corresponding  $3$-torus in the canonical order.
Now, define $Y^{}_{\varepsilon}$ as the quotient of the disjoint union of the tori and cylinders by the equivalence relation defined above.
Endowed with the induced length metric this construction yields, for any $\varepsilon \geqslant 0$, a compact connected non-positively curved geodesic space $Y^{}_{\varepsilon}$.
Its universal cover $X^{}_{\varepsilon} = \widetilde{Y}^{}_{\varepsilon}$ is a locally compact Hadamard space.

By stretching the interval factor of the cylinders $C^{\pm}_2$ in $Y^{}_0$ by $(1 + 2\varepsilon)$ and mapping the $3$-tori accordingly we obtain a mapping from $Y^{}_0$ into $Y^{}_{\varepsilon}$.
This mapping is one to one on all the vertices but $V^{\pm}_2$.
The $3$-tori therein are mapped into the $3$-torus-fibers over $0$ in the interval factor of the thickened $3$-tori in $Y^{}_{\varepsilon}$.
The lifting of this homotopy equivalence to the universal covers is a quasi-isometry which we denote by $\Phi^{}_{\varepsilon}$.

\subsection{Admissibility and geometric data}
\label{sec:cex:adm_geom_data}
In order to show that the action of $\Gamma = \pi^{}_1 (Y^{}_{\varepsilon})$ is admissible we consider $Y^{}_{\varepsilon}$ as a graph of spaces with eight vertices $V^{\pm}_0 = E^{\pm}_0 \cup C^{\pm}_0 \cup E^{\pm}_1 \cup C^{\pm}_{\times} \cup E^{}_{\times}$ and $V^{\pm}_i = E^{\pm}_i \cup C^{\pm}_i \cup E^{\pm}_{i+1}$ for $i = 1,2,3$ and nine edges $E^{}_{\times} = V^-_0 \cap V^+_0$ and $E^{\pm}_i = V^{\pm}_{i-1} \cap V^{\pm}_i$ for $i = 0,1,2,3$.
Then $\Gamma$ has a graph of groups decomposition over the graph $\G$
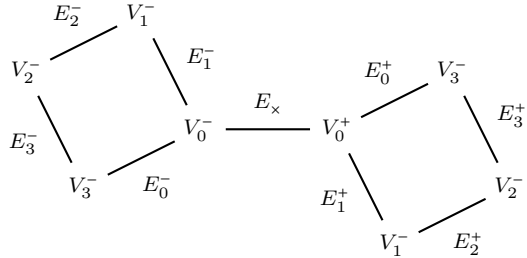
\begin{equation} \label{eq:Yd_graph_of_spaces}
\begin{pspicture}(-4,-1.7)(4,1.7)
	\psset{xunit=2em, yunit=2em}
	\cnodeput*{0}(-2.2,2){a1}{$\scriptstyle V^-_1$}
	\cnodeput*{0}(-4.2,1){a2}{$\scriptstyle V^-_2$}
	\cnodeput*{0}(-3.2,-1){a3}{$\scriptstyle V^-_3$}
	\cnodeput*{0}(2.2,-2){b1}{$\scriptstyle V^-_1$}
	\cnodeput*{0}(4.2,-1){b2}{$\scriptstyle V^-_2$}
	\cnodeput*{0}(3.2,1){b3}{$\scriptstyle V^-_3$}
	\cnodeput*{0}(-1.2,0){a0}{$\scriptstyle V^-_0$}
	\cnodeput*{0}(1.2,0){b0}{$\scriptstyle V^+_0$}
	\ncline{-}{a0}{b0} \naput{$\scriptstyle E^{}_{\times}$}
	\ncline{-}{a0}{a1} \nbput{$\scriptstyle E^-_1$}
	\ncline{-}{a1}{a2} \nbput{$\scriptstyle E^-_2$}
	\ncline{-}{a2}{a3} \nbput{$\scriptstyle E^-_3$}
	\ncline{-}{a3}{a0} \nbput{$\scriptstyle E^-_0$}
	\ncline{-}{b0}{b1} \nbput{$\scriptstyle E^+_1$}
	\ncline{-}{b1}{b2} \nbput{$\scriptstyle E^+_2$}
	\ncline{-}{b2}{b3} \nbput{$\scriptstyle E^+_3$}
	\ncline{-}{b3}{b0} \nbput{$\scriptstyle E^+_0$}
\end{pspicture}
\end{equation}
with egde- and vertex-groups $\Gamma^{}_{\gr{G}} (E^{}_{\times}) \cong \Gamma^{}_{\gr{G}} (E^{\pm}_i) = \pi^{}_1 (E^{\pm}_i) \cong \ZZ^3$, $\Gamma^{}_{\gr{G}} (V^{\pm}_i) = \pi^{}_1 (V^{\pm}_i) \cong \FF_2 \times \ZZ^2$ and $\Gamma^{}_{\gr{G}} (V^{\pm}_0) = \pi^{}_1 (V^{\pm}_0) \cong \FF_3 \times \ZZ^2$.
The center of each vertex group is generated by the fundamental group of the cylinder%
\footnote{Although the vertices $V^-_0$ and $V^+_0$ each contain two cylinders, their respective fundamental groups coincide in $\pi^{}_1 (V^{\pm}_0)$.}
along which the $3$-tori are glued, so define $\Lambda^{}_{\gr{G}} (V^{\pm}_0) = \pi^{}_1 (C^{\pm}_0) \cong \ZZ^2$.
Any edge group is generated by the centers of its adjacent vertex groups and its image in either vertex group is generated by its center and a generator of its free factor.
For any two adjacent edges these elements of the free factor do not coincide.
Therefore, the images of any two groups carried by consecutive edges are not commensurable.
Every quotient $\rquotient{\pi^{}_1 (V^{\pm}_i)}{\pi^{}_1 (C^{\pm}_i)}$ is a free group on either two or three generators; hence it is non-elementary hyperbolic.
In summary, we have shown:

\begin{prop} \label{cex:prop:adm}
The action of $\Gamma$ on $X^{}_{\varepsilon}$ is admissible of rank $3$.
\end{prop}

If $v$ is a lift of a vertex $V^{\pm}_i$ to the Bass-Serre tree of $\Gamma_{\G} (\_)$, the free abelian subgroup $\Lambda^{}_v$ is the lift of the center of the vertex group into the stabilizer $\Gamma^{}_v$.
Its minimal displacement set $\minset_{X^{}_{\varepsilon}} (\Lambda^{}_v)$ is isometric to the universal cover of $V^{\pm}_i$.
Therefore we get a $\Gamma^{}_v$-equivariant isometric embedding of $\minset_{X^{}_0} (\Lambda^{}_v)$  into $\minset_{X^{}_{\varepsilon}} (\Lambda^{}_v)$ by identifying the $3$-tori of the vertex spaces $V^{\pm}_2 \subset Y^{}_0$ with the innermost $3$-tori of the thickened vertex spaces $V^{\pm}_2 \subset Y^{}_{\varepsilon}$.
Since the generalized geometric data are determined by the actions of the vertex stabilizers on these minimal displacement sets, we have proved

\begin{prop} \label{prop:X0_Xd_geom_data}
For any $\varepsilon > 0$ the generalized geometric data of the $\Gamma$-actions on $X^{}_0$ and $X^{}_{\varepsilon}$ are equivalent.
\end{prop}

\subsection{Geometric decomposition}
\label{sec:cex:geom_decomp}
Let $\T$ denote the Bass-Serre tree of $\Gamma_{\G} (\_)$.
The space $X^{}_{\varepsilon}$ is the union of walls $\W^{}_e$, one for each edge $e \in \T$, and strips $\S_{(e , e')}$ for consecutive edges $e , e' \in \T$:
For any edge $e \in \T$ the wall $\W^{}_e$ is a copy of the universal cover of $E = \Gamma . e \in \E\G$, which is isometric to either $\EE^3$ or $\EE^3 \times [0,\varepsilon]$, and we call it of type $E$.
Each strip $\S_{(e , e')}$ is a copy of the universal cover $\EE^2 \times [0,1]$ of the cylinder $C$ contained in the vertex $V = \Gamma . (e \cap e') \in \V\G$ and shall be called of type $C$.
How each strip is glued along its boundary components to exactly two walls can be seen as follows:
We identify $\G$ with a fundamental domain of the $\Gamma$-action on $\T$ and consider $\pi^{}_1 (V) = \Gamma^{}_{\G} (V)$ and $\pi^{}_1 (E) = \Gamma^{}_{\G} (E)$ to be subgroups of $\Gamma$.
Now, let $e$ and $e'$ be edges which a common vertex $v$.
Then there are $\gamma , \gamma' \in \Gamma$,  $E , E' \in \E\G$ and $V \in \V\G$ such that $e = \gamma . E$, $v = \gamma . V$, $e' = \gamma' . E'$ and $\gamma' \gamma^{-1} \in \Gamma^{}_V$; the elements $\gamma$ and $\gamma'$ are unique up to $\Lambda^{}_e \cap \Lambda^{}_{e'} = \Lambda^{}_v$.
Since $E$ and $E'$ are adjacent, they are incident to a cylinder $C \subset V$.

Each wall $\W^{}_e$ is the minimal displacement set of $\Lambda_e$ which acts cocompactly by translations on the $\EE^3$-factor and trivially on the interval-factor in case $e$ is a thickened wall, i.e., of type $E^{\pm}_2$ or $E^{\pm}_3$.
The boundary component of $\S_{(e , e')}$ covering $E \cap C$ is precisely the convex hull of a $\Lambda^{}_v$-orbit in $\W_e$.
We refer to this two-dimensional affine subspace as the fringe $\L$ of $\S_{(e , e')}$ in $\W^{}_e$.
Let $e_1$, $e_2$ and $e_3$ be subsequent edges in $\T$, $v_1$ and $v_2$ the vertices in between and $\S_i = \S_{(e_i , e_{i+1})}$ the strips incident to the wall $\W = \W_{e_2}$.
Since, by Proposition \ref{cex:prop:adm}, $\smash{\Lambda_{v_1}}$ and $\smash{\Lambda_{v_2}}$ are incommensurable subgroups of $\smash{\Lambda_{e_2}}$, the fringes $\smash{\L^-_{\W} = \S_1 \cap \W}$ and $\smash{\L^+_{\W} = \S_2 \cap \W}$ are, possibly after an appropriate projection onto the $\EE^3$-factor of $\W$, non-parallel affine hyperplanes; hence they intersect in a line.
In case $\W \cong \EE^3$, define $\smash{\L^{\cap}_{\W} = \L^-_{\W} \cap \L^+_{\W}}$, and if $\W$ is of either type $E^{\pm}_2$ or type $E^{\pm}_3$, denote the intersection of the appropriate projections of the fringes along the interval factor of the wall by $\smash{\L^{\cap -}_{\W} \hookrightarrow \EE^3 \times \{ 0\}}$ and $\smash{\L^{\cap +}_{\W} \hookrightarrow \EE^3 \times \{ \varepsilon \}}$.

The convex subspace $Y^-_{\varepsilon} = V^-_0 \cup V^-_1 \cup V^-_2 \cup V^-_3$ in $Y^{}_{\varepsilon}$ is a metric product $\check{Y}^-_{\varepsilon} \times \sphere^1$:
The $\gamma^{}_i$ and $\gamma^-_{\times}$ generate an $\sphere^1$-factor in each torus or cylinder in $Y^-_{\varepsilon}$, and according to the relations \eqref{eq:Yd_gluing_koord} these spheres are parallel.
Similar to $Y^{}_{\varepsilon}$ itself, its subspace $\check{Y}^-_{\varepsilon}$ consists of five (thickened) flat $2$-tori $\check{E}^-_i$ and $\check{E}^{}_{\times}$ which are amalgamated along the boundaries of five cylinders $\check{C}^-_i$ and $\check{C}^-_{\times}$.
This defines a subgraph $\G^-$ of $\G$ with four vertices $\check{V}^-_0 = \check{E}^-_0 \cup \check{C}^-_0 \cup \check{E}^-_1 \cup \check{C}^-_{\times} \cup \check{E}^{}_{\times}$, $\check{V}^-_i  = \check{E}^-_i \cup \check{C}^-_i \cup \check{E}^-_{i+1}$ and four edges 
$\check{E}^-_i = \check{V}^-_{i-1} \cap \check{V}^-_i$.
Let $\check{\Gamma}_{\G^-} (\_)$ denote the associated graph of groups.

\begin{rem*}
Anticipating a notion introduced in section \ref{sec:salvage} we have $\Lambda^{\cap}_e = \langle \gamma^-_{\times} \rangle$ for a lift $e$ of $E^-_1$.
Therefore the subtree $\check{\T} (e)$, to be defined in Proposition \ref{prop:reduced_adm_action}, is isomorphic to the Bass-Serre tree of $\check{\Gamma}^{}_{\G^-} (\_)$.
\end{rem*}

The canonical generators of the fundamental groups of the $2$-tori, represented by the coordinate lines, will still be denoted by $\alpha^-$ and $\beta^-$, indexed according to the edge containing them.
For the fundamental group $\check{\Gamma}$ and the universal cover $\check{X}^-_{\varepsilon}$ of $\check{Y}^-_{\varepsilon}$ we have that $\check{X}^-_{\varepsilon} \times \EE \cong \minset_{X^{}_{\varepsilon}} (\gamma^-_{\times})$ and $\check{\Gamma} \cong \pi^{}_1 \left( \check{\Gamma}^{}_{\G^-} (\_) \right) \cong \rquotient{\centralizer_{\Gamma} (\gamma^-_{\times})}{\langle \gamma^-_{\times} \rangle}$.
Hence $\check{X}^-_{\varepsilon}$ is a locally compact Hadamard space upon which $\check{\Gamma}$ acts admissibly of rank $2$.
For any $\varepsilon \geqslant 0$ the geometric data on $\check{X}^-_0$ and $\check{X}^-_{\varepsilon}$ are equivalent.
In fact, we may regard $\check{X}^-_{\varepsilon}$ as a convex subspace of $X^{}_{\varepsilon}$.
Any two lifts of $\check{X}_{\varepsilon}$ to $X_{\varepsilon}$ are isometric by either a multiple of $\gamma^+_{\times}$ or an element in $\Gamma \setminus \centralizer_{\Gamma} (\gamma^+_{\times})$.

Analogous to the geometric decomposition we discussed at the beginning of this section, we decompose $\check{X}^-_{\varepsilon}$ into walls $\check{\W} \cong \EE^2$ or $\check{\W} \cong \EE^2 \times [0,\varepsilon]$ and strips $\check{\S} \cong \EE \times [0,1]$.
Let $\check{\L}^-_i$ and $\check{\L}^+_i$ the fringes of two strips $\check{\S}^{}_{i-1}$ und $\check{\S}^{}_i$ incident to a wall $\check{\W}^{}_i$.
Denote the intersection of their projections to the $\EE^2$-factor of the wall by $\check{o}^{}_i$ or $\check{o}^-_i \in \check{\L}^-_i$ and $\check{o}^+_i \in \check{\L}^+_i$, depending on whether $\check{\W}^{}_i$ has a trivial intervall factor or not.
Since the fringes $\check{\L}^-_i$ are axes of the central elements in the vertex-group, they are oriented by our choice \eqref{eq:Yd_gluing_koord} of representatives in the cylinders:
A point in a fringe is said to lie above $\check{o}^{}_i$ if it is a positive multiple of the translation by the chosen generator of the center.

The spaces $Y^+_{\varepsilon}$, $\check{Y}^+_{\varepsilon}$, $\check{X}^+_{\varepsilon}$ and the subgraph $\G^+$ are defined accordingly and we have $\check{\Gamma} \cong \pi^{}_1 (\check{Y}^-_{\varepsilon}) \cong \pi^{}_1 (\check{Y}^+_{\varepsilon})$.

The quasi-isometry $\Phi_{\varepsilon}$ obviously preserves the type of each wall.
In particular, for thickened walls $\Phi_{\varepsilon}$ coincides with the canonical embedding $\EE^3 \times \{ 0 \} \hookrightarrow \EE^3 \times [0,\varepsilon]$.
Moreover, $\Phi^{}_{\varepsilon}$ is compatible with the decomposition $\check{X}^{\pm}_0 \times \EE$ and splits as a product $\check{\Phi}^{\pm}_{\varepsilon} \times \id^{}_{\EE}$; its first factor is a $\check{\Gamma}$-equivariant quasi-isometry which, for sake of simplicity, we will denote by $\Phi^{}_{\varepsilon}$ as well.

\begin{rem*}
The quasi-isometry $\Phi_{\varepsilon}$ is a universal template map in the following sense:
Since each wall in $X_{\varepsilon}$ is two-sided, we may associate to any geodesic ray the sequence of types of walls it intersects.
The subspace of the corresponding walls and strips in $X_0$ is a template in the sense of \cite[Chapter 4]{crkl_geod} with the restriction of $\Phi^{}_{\varepsilon}$ constituting its template map.
\end{rem*}

\subsection{Special segments}
\label{sec:cex:spec_segm}
We will now construct a sequence of segments which shall later be utilized in proving that, for positive $\delta$, the quasi-isometry $\Phi^{}_{\delta}$ does not extend continuously to the virtual boundary of $X_0$.

\begin{figure}[!ht]
\centering
\begin{pspicture}(-5.882,-2.508)(2.299,6.494)
\psset{linejoin=1}
\psline[linewidth=0.8pt,linestyle=dashed,dash=1pt 1pt,arrowsize=3pt 0,ArrowInside=->](0,0)(2.09,0)
\psline[linewidth=0.3pt,arrowsize=3pt 0,ArrowInside=->](2.09,0)(2.09,2.09)
\psline[linewidth=0.8pt,linestyle=dashed,dash=1pt 1pt,arrowsize=3pt 0,ArrowInside=->](0,2.09)(2.09,2.09)
\psline[linewidth=0.3pt,arrowsize=3pt 0,ArrowInside=->](0,0)(0,2.09)
\psline[linewidth=0.6pt,linestyle=dashed,dash=4pt 2pt,arrowsize=3pt 0,ArrowInside=->](0,2.09)(2.09,0)
\psdots[dotstyle=*,dotsize=4pt](0,2.09)
\psline[linewidth=0.8pt,linestyle=dashed,dash=1pt 1pt,arrowsize=3pt 0,ArrowInside=->](1.478,5.016)(1.478,6.494)
\psline[linewidth=0.8pt,linestyle=dashed,dash=1pt 1pt,arrowsize=3pt 0,ArrowInside=->](0,5.016)(0,6.494)
\psline[linewidth=0.6pt,linestyle=dashed,dash=4pt 2pt,arrowsize=3pt 0,ArrowInside=->](0,5.016)(1.478,6.494)
\psline[linewidth=0.3pt,arrowsize=3pt 0,ArrowInside=->](0,6.494)(1.478,6.494)
\psline[linewidth=1.2pt,arrowsize=3pt 2,ArrowInside=->,ArrowInsideNo=2](0,5.016)(1.478,5.016)
\psdots[dotstyle=*,dotsize=4pt](0,5.016)
\psdots[dotstyle=o,dotsize=4pt](1.478,5.016)
\psline[linewidth=0.6pt,linestyle=dashed,dash=4pt 2pt,arrowsize=3pt 0,ArrowInside=->](-2.926,5.016)(-2.926,6.494)
\psline[linewidth=0.6pt,linestyle=dashed,dash=4pt 2pt,arrowsize=3pt 0,ArrowInside=->](-4.404,5.016)(-4.404,6.494)
\psline[linewidth=0.8pt,linestyle=dashed,dash=1pt 1pt,arrowsize=3pt 0,ArrowInside=->](-4.404,5.016)(-2.926,6.494)
\psline[linewidth=0.3pt,arrowsize=3pt 0,ArrowInside=->](-4.404,6.494)(-2.926,6.494)
\psline[linewidth=1.2pt,arrowsize=3pt 2,ArrowInside=->,ArrowInsideNo=2](-2.926,5.016)(-4.404,5.016)
\psdots[dotstyle=*,dotsize=4pt](-2.926,5.016)
\psdots[dotstyle=o,dotsize=4pt](-4.404,5.016)
\psline[linewidth=0.6pt,linestyle=dashed,dash=4pt 2pt,arrowsize=3pt 0,ArrowInside=->](-5.016,0)(-2.926,0)
\psline[linewidth=0.3pt,arrowsize=3pt 0,ArrowInside=->](-2.926,0)(-2.926,2.09)
\psline[linewidth=0.3pt,arrowsize=3pt 0,ArrowInside=->](-5.016,0)(-5.016,2.09)
\psline[linewidth=0.8pt,linestyle=dashed,dash=1pt 1pt,arrowsize=3pt 0,ArrowInside=->](-5.016,2.09)(-2.926,0)
\psline[linewidth=0.6pt,linestyle=dashed,dash=4pt 2pt,arrowsize=3pt 0,ArrowInside=->](-5.016,2.09)(-2.926,2.09)
\psdots[dotstyle=*,dotsize=4pt](-5.016,2.09)
\psline[linewidth=0.8pt,linestyle=dashed,dash=1pt 1pt,arrowsize=3pt 0,ArrowInside=->](0,2.508)(2.09,2.508)
\psline[linewidth=0.3pt,arrowsize=3pt 0,ArrowInside=->](2.09,2.508)(2.09,4.598)
\psline[linewidth=0.6pt,linestyle=dashed,dash=4pt 2pt,arrowsize=3pt 0,ArrowInside=->](0,4.598)(2.09,4.598)
\psline[linewidth=0.3pt,arrowsize=3pt 0,ArrowInside=->](0,2.508)(0,4.598)
\psline[linewidth=1.2pt,arrowsize=3pt 2,ArrowInside=->,ArrowInsideNo=2](0,2.508)(2.09,4.598)
\psdots[dotstyle=o,dotsize=4pt](0,2.508)
\psdots[dotstyle=*,dotsize=4pt](2.09,4.598)
\psline[linewidth=0.6pt,linestyle=dashed,dash=4pt 2pt,arrowsize=3pt 0,ArrowInside=->](-5.016,2.508)(-2.926,2.508)
\psline[linewidth=0.3pt,arrowsize=3pt 0,ArrowInside=->](-2.926,4.598)(-2.926,2.508)
\psline[linewidth=0.8pt,linestyle=dashed,dash=1pt 1pt,arrowsize=3pt 0,ArrowInside=->](-5.016,4.598)(-2.926,4.598)
\psline[linewidth=0.3pt,arrowsize=3pt 0,ArrowInside=->](-5.016,4.598)(-5.016,2.508)
\psline[linewidth=1.2pt,arrowsize=3pt 2,ArrowInside=->,ArrowInsideNo=2](-2.926,4.598)(-5.016,2.508)
\psdots[dotstyle=o,dotsize=4pt](-2.926,4.598)
\psdots[dotstyle=*,dotsize=4pt](-5.016,2.508)
\psline[linewidth=0.8pt,linestyle=dashed,dash=1pt 1pt,arrowsize=3pt 0,ArrowInside=->](-.418,5.016)(-.418,6.494)
\psline[linewidth=0.3pt,arrowsize=3pt 0,ArrowInside=->](-.418,6.494)(-2.508,6.494)
\psline[linewidth=0.6pt,linestyle=dashed,dash=4pt 2pt,arrowsize=3pt 0,ArrowInside=->](-2.508,5.016)(-2.508,6.494)
\psline[linewidth=1.2pt,arrowsize=3pt 2,ArrowInside=->,ArrowInsideNo=2](-.418,5.016)(-2.508,5.016)
\psdots[dotstyle=o,dotsize=4pt](-.418,5.016)
\psdots[dotstyle=*,dotsize=4pt](-2.508,5.016)
\psline[linewidth=0.3pt,arrowsize=3pt 0,ArrowInside=->](-2.508,-.866)(-.418,-.866)
\psline[linewidth=0.6pt,linestyle=dashed,dash=4pt 2pt,arrowsize=3pt 0,ArrowInside=->](-.418,2.09)(-.418,-.866)
\psline[linewidth=0.8pt,linestyle=dashed,dash=1pt 1pt,arrowsize=3pt 0,ArrowInside=->](-2.508,2.09)(-2.508,-.866)
\psline[linewidth=1.2pt,arrowsize=3pt 2,ArrowInside=->,ArrowInsideNo=2](-2.508,2.09)(-.418,2.09)
\psdots[dotstyle=*,dotsize=4pt](-.418,2.09)
\psdots[dotstyle=o,dotsize=4pt](-2.508,2.09)
\psline[linewidth=0.6pt,linestyle=dashed,dash=4pt 2pt,arrowsize=3pt 0,ArrowInside=->](-5.882,-2.508)(-2.926,-2.508)
\psline[linewidth=0.3pt,arrowsize=3pt 0,ArrowInside=->](-2.926,-.418)(-2.926,-2.508)
\psline[linewidth=0.8pt,linestyle=dashed,dash=1pt 1pt,arrowsize=3pt 0,ArrowInside=->](-5.882,-.418)(-2.926,-.418)
\psline[linewidth=0.3pt,arrowsize=3pt 0,ArrowInside=->](-5.882,-.418)(-5.882,-2.508)
\uput[ul](0,2.09){$\scriptstyle \check{p}^{}_1$}
\uput[dl](2.09,-.209){$\scriptstyle \check{E}^-_1$}
\uput[d](1.045,0){$\scriptstyle \alpha^-_1$}
\uput[r](0,.836){$\scriptstyle \beta^-_1$}
\uput[ur](1.045,1.045){$\scriptstyle \alpha^-_1 {\beta^-_1}^{-1}$}
\uput[dl](0,5.016){$\scriptstyle \check{p}^-_2$}
\uput[u](.739,5.016){$\scriptstyle \hat{c}^-_0$}
\uput[ur](1.626,5.016){$\scriptstyle \check{E}^-_2$}
\uput[u](.739,6.494){$\scriptstyle \alpha^-_2$}
\uput[r](1.478,5.755){$\scriptstyle \beta^-_2$}
\uput[ul](.887,5.755){$\scriptstyle \alpha^-_2 \beta^-_2$}
\uput[dr](-2.926,5.016){$\scriptstyle \check{p}^+_3$}
\uput[u](-3.665,5.016){$\scriptstyle \hat{c}^-_0$}
\uput[ul](-4.552,5.016){$\scriptstyle \check{E}^-_3$}
\uput[u](-3.665,6.494){$\scriptstyle \alpha^-_3$}
\uput[l](-4.404,5.755){$\scriptstyle \beta^-_3$}
\uput[ul](-3.517,5.755){$\scriptstyle \alpha^-_3 \beta^-_3$}
\uput[ul](-5.016,2.09){$\scriptstyle \check{p}^{}_0$}
\uput[ul](-5.225,0){$\scriptstyle \check{E}^-_0$}
\uput[u](-3.971,0){$\scriptstyle \alpha^-_0$}
\uput[l](-5.016,1.045){$\scriptstyle \beta^-_0$}
\uput[ur](-3.971,1.045){$\scriptstyle \alpha^-_0 {\beta^-_0}^{-1}$}
\uput[dr](1.045,3.553){$\scriptstyle \hat{c}^-_0$}
\uput[ur](2.299,2.508){$\scriptstyle \check{C}^-_1$}
\uput[dr](-3.971,3.553){$\scriptstyle \hat{c}^-_0$}
\uput[ul](-5.225,2.508){$\scriptstyle \check{C}^-_3$}
\uput[u](-1.463,5.016){$\scriptstyle \hat{c}^-_0$}
\uput[dr](-2.508,4.807){$\scriptstyle \check{C}^-_2$}
\uput[d](-1.463,2.09){$\scriptstyle \hat{c}^-_0$}
\uput[dl](-.418,-1.075){$\scriptstyle \check{C}^-_0$}
\uput[ur](-2.717,-2.508){$\scriptstyle \check{C}^-_{\times}$}
\end{pspicture}
\caption[$\hat{c}^-_0$ in the fundamental domains in $\check{Y}^-_0$]{The dashed or dotted lines of adjacent fundamental domains are identified. $\smash{\hat{c}^-_0}$ is drawn bold and decorated with double arrows.}
\label{fig:Y0-_domains}
\end{figure}
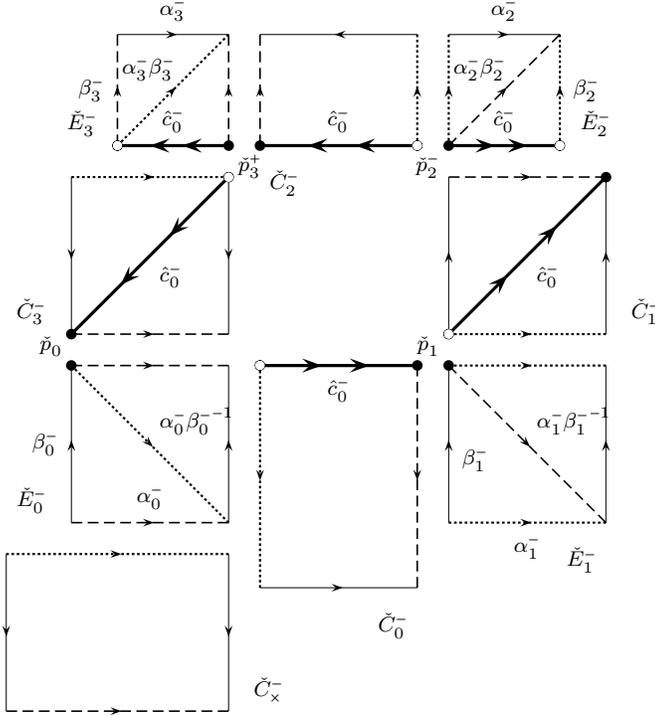

First of all we define a geodesic ray in $\check{X}^-_{\varepsilon}$ whose projection to $\check{X}^-_{\varepsilon}$ intersects the edges $\check{E}^-_i$ in cyclic order.
To that end let $p^-_2 = (0,0,0;0)^-_2 \in E^-_2$, $p^+_2  = (0,0,0;\varepsilon)^-_2 \in E^-_2$, $p^-_3 = (0,0,0;\varepsilon)^-_3 \in E^-_3$, $p^+_3  = (0,0,0;0)^-_3 \in E^-_3$, and $p^{}_i = (0,0,0;0)^-_i \in E^-_i$ for $i \in \{ 0,1 \}$ and denote the respective projections to $\check{E}^{\pm}_i$ by $\check{p}^{\pm}_i$; then $\overline{\check{p}^-_2 \check{p}^+_2} \cup \overline{\check{p}^+_2 \check{p}^-_3} \cup \overline{\check{p}^-_3 \check{p}^+_3} = \overline{\check{p}^-_2 \check{p}^+_3}$.
We define $\hat{c}^-_{\varepsilon}$ to be the geodesic representative%
\footnote{This corresponds to $(\alpha^-_2 \beta^-_2)^2 \smash{\tau_{\check{E}^-_2}} (\alpha^-_3 \beta^-_3)^2$ in the presentation of $\check{\Gamma} \cong \pi^{}_1 ( \check{\Gamma}_{\G^-} (\_) ; \G^- \setminus \check{E}^-_2 )$.}
of
\begin{equation} \label{eq:c-_rep}
\smash{\left[ \overline{\check{p}^{}_0 \check{p}^{}_1} \cup \overline{\check{p}^{}_1 \check{p}^-_2} \cup (\alpha^-_2 \beta^-_2)^2 \cup \overline{\check{p}^-_2 \check{p}^+_3} \cup (\alpha^-_3 \beta^-_3)^{-2}_{\vphantom{0}} \cup \overline{\check{p}^+_3 \check{p}^{}_0} \right] \in \pi_1 (\check{Y}^-_{\varepsilon})} \; .
\end{equation}
Figure \ref{fig:Y0-_domains} depicts the loop $\hat{c}^-_0$ in the fundamental domains of relevant $2$-tori and cylinders in $\check{Y}^-_0$.
%

For a lift%
\footnote{Any lift of $\check{p}^{\pm}_i$ is an intersection of two fringes since $\check{p}^{\pm}_i$ is the intersection of two cylinders adjacent to $\check{E}^{\pm}_i$.}
$\check{o}^{}_0 \in \pi^{-1} (\check{p}^{}_0)$ we let $\check{c}^-_{\varepsilon}$ be a lift of $\hat{c}^-_{\varepsilon}$ to $\check{X}^-_{\varepsilon}$ starting at $\check{o}^{}_0$. 
Then the geodesic ray $\check{c}^-_{\varepsilon}$ intersects a sequence of walls $\check{\W}^{}_0 , \check{\W}^{}_1 , \check{\W}^{}_2 , \check{\W}^{}_3 , \check{\W}^{}_4 , \check{\W}^{}_5 , \ldots$ of type $\check{E}^-_0 , \check{E}^-_1 , \check{E}^-_2 , \check{E}^-_3 , \check{E}^-_0 , \check{E}^-_1 , \ldots $ resp.\ in order.
Denote the fringe between $\check{\W}^{}_i$ and $\check{\W}^{}_{i+1}$ by \nolinebreak $\check{\S}^{}_i$.

\begin{figure}[h!]
\centering
\makebox[\textwidth][c]{\input{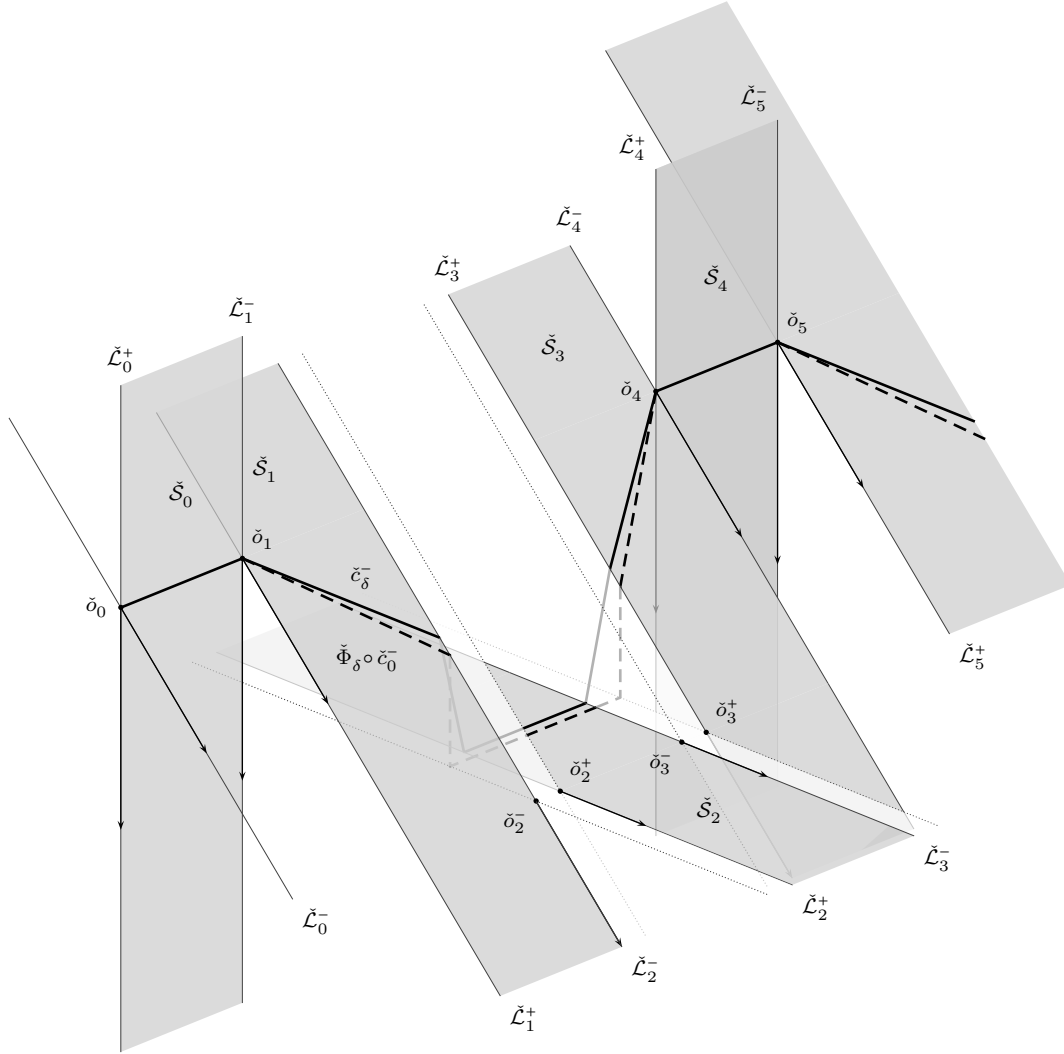}}
\caption{
The arrows indicate the translation vectors of the generators orienting the strips.
The intersections of the fringes in strips with odd indices are shifted along the real factor by twice the translation length of either $\alpha^-_2 \beta^-_2 \sim \alpha^-_1$ or $\alpha^-_3 \beta^-_3 \sim \alpha^-_0$.
}
\label{fig:example_segment_dim_3}
\end{figure}

In a next step we will compare the geodesic ray $\check{c}^-_{\delta}$ with the quasi-isometric image of $\check{c}^-_0$ in the universal cover $\smash{\check{X}^-_{\delta}}$.
We will show, in particular, that $\hat{c}^-_0$ is strictly shorter than $\hat{c}^-_{\delta}$.
Figure \ref{fig:example_segment_dim_3} is a draft of the walls and strips the ray $\check{c}^-_{\delta}$ traverses and the run of $\check{c}^-_{\delta}$ and $\Phi^{}_{\delta} \circ \check{c}^-_0$ therein.
%
%
Since $\check{X}^-_{\varepsilon}$ is piecewise Euclidean, any geodesic is linear on the walls and strips.
As is easily seen from the relations \ref{eq:Yd_gluing_koord}, the fringes in any wall intersect at an angle of $\pm \frac{\pi}{4}$.
For reasons of symmetry the geodesic segment from $\check{o}^{}_0$ to $\check{o}_5$ passes the center axis $\check{\temp{C}}$ of the strip $\check{\S}_2$ at an angle of $\frac{\pi}{2}$ (see Figures \ref{fig:example_segment_dim_3} and \ref{fig:o1_C}).
%
\begin{figure}[H]
\centering
\makebox[\textwidth][c]{
\begin{pspicture}(-1.042,-4.5)(6.114,1.899)
\psset{linejoin=1}
\pspolygon[linewidth=0,linestyle=none,fillcolor=lightgray!80,fillstyle=solid,opacity=0.7](.563,-1.412)(4.38,-3.383)(6.114,-2.841)(2.298,-.87)
\psline[linewidth=0.9pt,linestyle=dotted,dotsep=1pt](1.433,-1.142)(4.293,-2.619)
\pspolygon[fillcolor=white,fillstyle=solid,opacity=0.7,linestyle=none](3.081,-3.004)(3.558,-3.962)(3.905,-3.853)(3.428,-2.895)
\psline[linewidth=0,linestyle=none](3.081,-3.004)(3.558,-3.962)(3.905,-3.853)(3.428,-2.895)
\psline[linewidth=0.5pt,linecolor=darkgray,linestyle=dotted,dotsep=0.7pt](3.431,-2.897)(3.905,-3.849)
\pspolygon[linewidth=0,linestyle=none,fillcolor=white,fillstyle=solid,opacity=0.7](3.085,-3.003)(3.665,-4.168)(4.245,-3.602)
\pspolygon[linewidth=0,linestyle=none,fillcolor=white,fillstyle=solid,opacity=0.7](.219,-1.526)(4.035,-3.497)(4.382,-3.388)(.566,-1.417)
\psline[linewidth=0.3pt,linecolor=darkgray](.566,-1.413)(3.425,-2.89)
\psline[linewidth=0.6pt](2.013,-.858)(2.359,-2.34)(3.227,-2.069)
\psdots[dotstyle=*,dotsize=2pt](2.359,-2.34)
\psline[linewidth=0.5pt,linestyle=dashed](2.359,-.749)(2.359,-2.34)(2.012,-2.448)(2.012,-.858)
\pspolygon[linewidth=0,linestyle=none,fillcolor=white,fillstyle=solid,opacity=0.7](.219,-1.526)(3.081,-3.004)(.696,1.786)
\psline[linewidth=0.5pt,linecolor=darkgray,linestyle=dotted,dotsep=0.7pt](.219,-1.522)(3.078,-2.998)
\psdots[dotstyle=*,dotsize=2pt](3.227,-2.069)
\pspolygon[fillcolor=white,fillstyle=solid,opacity=0.7,linestyle=none](.696,1.786)(3.078,-2.998)(3.425,-2.89)(1.043,1.895)
\psline[linewidth=0,linestyle=none](3.425,-2.89)(1.043,1.895)(.696,1.786)(3.078,-2.998)
\psline[linewidth=0.5pt,linecolor=darkgray,linestyle=dotted,dotsep=0.7pt](1.043,1.899)(3.428,-2.891)
\psline[linewidth=0.4pt](1.52,.941)(1.693,.995)
\psline[linewidth=0.4pt,arrowscale=0.4,arrows=t-t](1.624,.973)(2.463,-.713)
\pspolygon[fillcolor=lightgray!80,fillstyle=solid,opacity=0.7,linestyle=none](-1.042,1.25)(1.343,-3.54)(3.078,-2.998)(.693,1.792)
\psline[linewidth=0,linestyle=none](3.078,-2.998)(.693,1.792)(-1.042,1.25)(1.343,-3.54)
\pspolygon[fillcolor=white,fillstyle=solid,opacity=0.7,linestyle=none](3.078,-2.998)(3.081,-3.004)(3.428,-2.895)(3.425,-2.89)
\psline[linewidth=0,linestyle=none](3.078,-2.998)(3.081,-3.004)
\psline[linewidth=0,linestyle=none](3.428,-2.895)(3.425,-2.89)
\psline[linewidth=0.3pt,linecolor=darkgray](.696,1.79)(3.081,-3)
\psline[linewidth=0.5pt,linestyle=dashed](1.173,.832)(1.52,.941)
\psline[linewidth=0.5pt,linestyle=dashed](-.562,.29)(1.173,.832)
\psline[linewidth=0.4pt](2.359,-.745)(2.533,-.691)
\psline[linewidth=0.9pt,linestyle=dotted,dotsep=1pt](4.293,-2.619)(5.25,-3.113)
\psline[linewidth=0.5pt,linestyle=dashed](2.012,-.854)(2.359,-.745)(2.359,-.749)
\psline[linewidth=0.6pt](2.012,-.854)(2.013,-.858)
\psline[linewidth=0.5pt,linestyle=dashed](2.012,-.858)(2.012,-.854)
\psline[linewidth=0.6pt](-.562,.29)(2.012,-.854)
\psdots[dotstyle=*,dotsize=2pt](2.012,-.854)
\psline[linewidth=0.3pt,linecolor=darkgray](3.425,-2.89)(4.382,-3.384)
\psline[linewidth=0.5pt,linecolor=darkgray,linestyle=dotted,dotsep=0.7pt](3.428,-2.891)(3.431,-2.897)
\psdots[dotstyle=*,dotsize=2pt](3.428,-2.891)
\pspolygon[fillcolor=lightgray!80,fillstyle=solid,opacity=0.7,linestyle=none](1.343,-3.54)(1.346,-3.546)(3.081,-3.004)(3.078,-2.998)
\psline[linewidth=0,linestyle=none](1.343,-3.54)(1.346,-3.546)
\psline[linewidth=0,linestyle=none](3.081,-3.004)(3.078,-2.998)
\pspolygon[fillcolor=lightgray!80,fillstyle=solid,opacity=0.7,linestyle=none](1.346,-3.546)(1.82,-4.498)(3.555,-3.956)(3.081,-3.004)
\psline[linewidth=0,linestyle=none](1.346,-3.546)(1.82,-4.498)(3.555,-3.956)(3.081,-3.004)
\psline[linewidth=0.3pt,linecolor=darkgray](-1.039,1.248)(1.349,-3.547)
\psline[linewidth=0.5pt,linecolor=darkgray,linestyle=dotted,dotsep=0.7pt](3.078,-2.998)(4.035,-3.492)
\psline[linewidth=0.3pt,linecolor=darkgray](3.081,-3)(3.084,-3.005)
\psdots[dotstyle=*,dotsize=2pt](3.081,-3)
\psline[linewidth=0.3pt,linecolor=darkgray](3.084,-3.005)(3.558,-3.958)
\psdots[dotstyle=*,dotsize=2pt](-.562,.29)
\psline[linewidth=0.3pt,linecolor=darkgray](1.349,-3.547)(1.823,-4.5)
\uput[ur](3.227,-2.069){$\scriptstyle \check{q}^{}_{\check{\temp{C}}} (s)$}
\uput[r](2.123,.133){$\scriptstyle s$}
\uput[dl](2.012,-.854){$\scriptstyle \check{q}^-(s)$}
\uput[l](-.562,.29){$\scriptstyle \check{o}^{}_1$}
\uput[l](3.081,-3){$\scriptstyle \check{o}^-_2$}
\uput[ur](3.428,-2.891){$\scriptstyle \check{o}^+_2$}
\uput[u](.028,.917){$\scriptstyle \check{\S}^{}_1$}
\uput[dl](5.64,-2.596){$\scriptstyle \check{\S}^{}_2$}
\uput[d](3.558,-3.958){$\scriptstyle \check{\L}^-_2$}
\uput[dr](4.382,-3.384){$\scriptstyle \check{\L}^+_2$}
\uput[dr](5.25,-3.113){$\scriptstyle \check{\temp{C}}$}
\end{pspicture}
}
\caption{\mbox{A broken geodesic parametrised by the distance $s$ that $\check{q}^- (s)$ lies above $\check{o}_1$.}}
\label{fig:o1_C}
\end{figure}
\removelastskip
\smallskip
First of all, consider the geodesic segment from $\check{o}^{}_1$ and the orthogonal projection $\smash{q_{\check{\temp{C}}} = \proj_{\check{\temp{C}}} (\check{o}^{}_1)}$ of $\check{o}^{}_1$ onto $\check{\temp{C}}$.
Any geodesic segment from $\check{o}_1$ to $\check{\temp{C}}$ will intersect the fringe $\smash{\check{\L}^-_2}$ at a point $\check{q}^-$.
Thus the geodesic segment from a point $\check{q}^- \in \check{\L}^-_2$ to its projection $\proj_{\check{\temp{C}}} (\check{q}^-)$ onto $\check{\temp{C}}$ traverses the wall $\check{\W}^{}_2$ and the strip $\check{\S}^{}_2$ orthogonally to the fringe $\check{\L}^+_2$.
It therefore suffices to minimize the length of the broken geodesic segment $\overline{\check{o}_1 \check{q}^- (s)} \cup \overline{\check{q}^- (s) \check{q}_{\check{\temp{C}}} (s)} \subset \check{X}^-_{\varepsilon}$ where $\check{q}^- (s)$ is the point on $\check{\L}^-_2$ lying at distance $s$ above the projection of $\check{o}_1$ onto $\check{\L}^-_2$ (see Figure \ref{fig:o1_C}).
A straightforward calculation shows the lengths $L (\varepsilon , s)$ of these paths to be strictly convex in $s$.
Therefore $L (\varepsilon , \, .)$ attains a minimum at some $s_{\varepsilon} \in (0,1]$.
Furthermore, $L (\varepsilon , s_{\varepsilon} )$ is strictly increasing in $\varepsilon$.

Secondly, an easy \textit{billiard}-argument shows that for any $s \in [-1 , 1]$ the segment $\overline{\check{o}_0 \check{q}^-(s)}$ does indeed contain the point $\check{o}_1$ (see Figure \ref{fig:billard}).
\medskip
\enlargethispage{2em}
\begin{figure}[H]
\centering
\begin{pspicture}(-2.3,-4.223)(3.081,3.138)
\psset{linejoin=1}
\pspolygon[linewidth=0,linestyle=none,fillcolor=lightgray!80,fillstyle=solid,opacity=0.7](-1.519,2.208)(1.343,-3.54)(3.078,-2.998)(.216,2.75)
\psline[linewidth=0.3pt,linecolor=darkgray](.219,2.748)(3.081,-3)
\psline[linewidth=0.4pt,linestyle=dashed](-1.23,-.583)(-.565,-1.918)
\psline[linewidth=0.4pt](1.364,.449)(1.537,.503)
\pspolygon[fillcolor=lightgray!80,fillstyle=solid,opacity=0.7,linestyle=none](-2.3,-.25)(-2.3,-3.097)(-.565,-2.555)(-.565,.292)
\psline[linewidth=0,linestyle=none](-2.3,-.25)(-2.3,-3.097)(-.565,-2.555)(-.565,.292)
\psline[linewidth=0.3pt,linecolor=darkgray](-.562,.286)(-.562,-2.556)
\psline[linewidth=0.4pt,linestyle=dashed](-.562,.29)(1.173,.832)
\psdots[dotstyle=*,dotsize=2pt](1.173,.832)
\psline[linewidth=0.4pt,arrowscale=0.4,arrows=t-t](1.468,.482)(2.708,-2.009)
\psline[linewidth=0.4pt,arrowscale=0.4,arrows=t-t](1.069,.8)(2.023,-1.116)
\psline[linewidth=0.4pt,linestyle=dashed](-.371,-.093)(1.364,.449)
\psline[linewidth=0.4pt,linestyle=dashed](-.565,-1.918)(.583,-4.223)
\psdots[dotstyle=*,dotsize=2pt](-.562,-1.924)
\psline[linewidth=0.6pt](-2.297,-.252)(-.562,-1.924)(.382,-3.962)
\psline[linewidth=0.4pt](1.953,-1.138)(2.127,-1.084)
\psline[linewidth=0.4pt](2.604,-2.042)(2.777,-1.988)
\psline[linewidth=0.3pt,linecolor=darkgray](-1.516,2.206)(-.565,.296)
\psline[linewidth=0.6pt](-.374,-.091)(-.371,-.093)(2.604,-2.042)
\psdots[dotstyle=*,dotsize=2pt](2.604,-2.042)
\psline[linewidth=0.4pt,arrowscale=0.4,arrows=t-t](-.267,.382)(-.076,-.001)
\psline[linewidth=0.6pt](-2.269,-.247)(-.562,.039)(-.374,-.091)
\psdots[dotstyle=*,dotsize=2pt](-.562,.039)
\pspolygon[fillcolor=lightgray!80,fillstyle=solid,opacity=0.7,linestyle=none](-2.3,2.596)(-2.3,-.25)(-.565,.292)(-.565,3.138)
\psline[linewidth=0,linestyle=none](-.565,.292)(-.565,3.138)(-2.3,2.596)(-2.3,-.25)
\psline[linewidth=0.3pt,linecolor=darkgray](-2.297,2.595)(-2.297,-3.098)
\psline[linewidth=0.3pt,linecolor=darkgray](-.562,3.137)(-.562,.286)
\psline[linewidth=0.3pt,linecolor=darkgray](-.565,.296)(1.346,-3.542)
\psdots[dotstyle=*,dotsize=2pt](-.562,.29)
\psline[linewidth=0.4pt,linestyle=dashed](-2.297,-.252)(-.562,.29)
\psdots[dotstyle=*,dotsize=2pt](-.371,-.093)
\psdots[dotstyle=*,dotsize=0pt](1.346,-3.944)
\psline[linewidth=0.6pt](-2.297,-.252)(-2.269,-.247)
\psdots[dotstyle=*,dotsize=2pt](-2.297,-.252)
\rput(-2.297,-.252){\psarc[linewidth=0.4pt]{->}{1.6}{9.5}{17.5}}
\uput[r](-1.048,.102){$\scriptstyle \theta$}
\rput(2.604,-2.042){\psarc[linewidth=0.4pt]{->}{0.7}{117}{147}}
\rput(-2.297,-.252){\psarc[linewidth=0.4pt]{->}{0.7}{316}{17.5}}
\rput(-.562,-1.924){\psarc[linewidth=0.4pt]{->}{0.7}{270}{294.8}}
\uput[300](-.562,-2.556){$\scriptstyle \leqslant \frac{\pi}{4}$}
\uput[r](-1.777,-.492){$\scriptstyle \geqslant \frac{\pi}{4}$}
\uput[l](-2.297,-.252){$\scriptstyle \check{o}^{}_0$}
\uput[ul](-.562,.29){$\scriptstyle \check{o}^{}_1$}
\uput[u](-.562,3.137){$\scriptstyle \check{\L}^-_1$}
\uput[dr](1.346,-3.542){$\scriptstyle \check{\L}^+_1$}
\uput[dr](3.081,-3){$\scriptstyle \check{\L}^-_2$}
\uput[r](1.173,.832){$\scriptstyle \proj^{}_{\check{\L}^-_2} (\check{o}^{}_1)$}
\uput[r](-.148,.265){$\scriptstyle > 0$}
\uput[l](1.572,-.372){$\scriptstyle 1$}
\uput[dr](2.149,-.523){$\scriptstyle \cot (\frac{\pi}{4} - \theta) > 1$}
\uput[ul](2.296,-1.584){$\scriptstyle \frac{\pi}{4} - \theta$}
\uput[dr](2.691,-2.015){$\scriptstyle \check{q}^- (s)$}
\uput[u](-1.43,.825){$\scriptstyle \check{\S}^{}_0$}
\uput[u](-.155,1.082){$\scriptstyle \check{\S}^{}_1$}
\end{pspicture}
\caption{Geodesics starting at $\smash{\check{o}_0}$ either intersect $\smash{\check{\L}^-_2}$ too far beyond $\smash{\check{o}_1}$ \mbox{(for small $\theta$)} or eventually stay in $\smash{\check{\W}_1}$ (for $\abs{\theta} \geqslant \smash{\frac{\pi}{4}}$).
}
\label{fig:billard}
\end{figure}
\pagebreak

The preceding discussion shows that the difference of the lengths $\check{l}^{}_{\delta} = l (\hat{c}^-_{\delta}) = 1 + 2 L (\delta , s^{}_{\delta})$ and $\check{l}^{}_0 = l(\hat{c}^-_0) = 1 + 2 L (0 , s^{}_0)$ of $\overline{\check{o}^{}_0 \check{o}^{}_4}$ in $\check{X}^-_{\delta}$ and $\check{X}^-_0$ resp.\ is strictly positive.

Utilizing the properties of the geodesic rays $\check{c}_{\varepsilon}$ in $\check{X}_{\varepsilon}$ we are about to construct a sequence of geodesic segments $c_{\varepsilon , n}$ in $X_{\varepsilon}$.
These segments start in $X^-_{\varepsilon}$ and end in $X^+_{\varepsilon}$.
In each subspace $c_{\varepsilon , n}$ will $n$-times intersect $E^-_0 , E^-_1 , E^-_2$ and $E^-_3$ in cyclic order.
Let $x^{}_{\varepsilon} = ( \check{c}^-_{\varepsilon} (0) , 0 ) \in \check{X}^-_{\varepsilon} \times \EE \subset X^{}_{\varepsilon}$ for $\varepsilon \geqslant 0$.
We may assume that $\Phi^{}_{\varepsilon} (x^{}_0) = x^{}_{\varepsilon}$.
For any $n \check{l}^{}_{\varepsilon}$ the geodesic segment $\smash{\check{c}^-_{\varepsilon}\!\!\mid_{[0,n \check{l}^{}_{\varepsilon}]}}$ ends at a point $\check{o}^-_{4n} \in \pi^{-1} (\check{p}^-_0)$, i.e., in a fringe of the strip $\check{\S}^-_{\times , n}$ of type $\check{C}^-_{\times}$.
It therefore can be extended, as a geodesic segment, by a segment of length $1$ into the strip $\smash{\check{\S}^-_{\times , n}}$ traversing it parallel to its interval factor.
Its endpoint $\check{q}^-_{\varepsilon , n}$ then lies in a fringe $\smash{\check{\L}^-_{\times , n}}$ in a wall $\smash{\check{\W}^{}_{\times , n}}$ of type $\smash{\check{E}^{}_{\times}}$.

The product path $\smash{\overline{\check{c}^-_{\varepsilon} (0) \check{q}^-_{\varepsilon , n}} \dtimes \id_{\EE}}$ is a geodesic segment of constant velocity and hits a fringe $\L^-_{\times , n}$ in a wall $\W^{}_{\times , n}$ of type $E^{}_{\times}$ at a point $q^-_{\varepsilon , n}$ at time $n \check{l}^{}_{\varepsilon} + 1$.
This may then be extended into the wall $\W^{}_{\times , n}$ traversing it orthogonally to the axes of $\alpha^{}_{\times}$ and at an angle of $\frac{\pi}{4}$ with both $\gamma^+_{\times}$ and $\gamma^-_{\times}$ (see Figure \ref{fig:Wx}).
Since $\Phi_{\varepsilon}$ maps $\W_{\times , n}$ isometrically to a wall of the same type, the quasi-isometric image of this geodesic segment is linear within $\W_{\times , n}$ and its angles with the axes are preserved.
The points $q^-_{\delta , n}$ and $\Phi^{}_{\delta} (q^-_{0 , n})$ lie on the same axis of $\gamma^-_{\times}$ in $\L^-_{\times , n}$, and each fringe covering $E^{}_{\times} \cap C^+_{\times}$ is orthogonal to the axis of $\gamma^-_{\times}$.
%
\begin{figure}[h!]
\centering
\begin{pspicture}(-4.952,-2.969)(3.095,5.469)
\psset{linejoin=1}
\pspolygon[linewidth=0,linestyle=none,fillcolor=lightgray!80,fillstyle=solid,opacity=0.7](-4.952,-.402)(-.619,-.939)(.928,1.408)(-3.404,1.944)
\pspolygon[linewidth=0,linestyle=none,fillcolor=lightgray!30,fillstyle=solid,opacity=0.7](-.619,-.939)(-.619,-2.969)(.928,-.623)(.928,1.408)
\psline[linewidth=0.5pt,linecolor=darkgray,linestyle=dotted,dotsep=0.7pt](-4.333,.536)(-.007,.001)
\psline[linewidth=1pt,linestyle=dashed](-3.25,.402)(-.007,3.04)
\psline[linewidth=1pt](-1.083,.134)(-.007,1.01)
\pspolygon[linewidth=0,linestyle=none,fillcolor=lightgray!80,fillstyle=solid,opacity=0.7](-.619,-.939)(.928,1.408)(.928,5.469)(-.619,3.123)
\psline[linewidth=0.5pt,arrows=->](0,-.011)(0,-2.031)
\pspolygon[linewidth=0,linestyle=none,fillcolor=lightgray!30,fillstyle=solid,opacity=0.7](-.619,-.939)(1.547,-1.207)(3.095,1.14)(.928,1.408)
\psline[linewidth=0.3pt,linecolor=darkgray](-.619,-.939)(.928,1.408)
\psline[linewidth=0.5pt,arrows=->](0,0)(.619,.939)
\psdots[dotstyle=*,dotsize=2pt](-3.25,.402)
\psdots[dotstyle=*,dotsize=2pt](-1.083,.134)
\psline[linewidth=0.5pt,arrows=-](0,0)(0,-.011)
\psline[linewidth=0.5pt,linecolor=darkgray,linestyle=dotted,dotsep=0.7pt](-.007,.001)(0,0)
\psline[linewidth=0.5pt,arrows=->](0,0)(2.166,-.268)
\psline[linewidth=0.5pt,linecolor=darkgray,linestyle=dotted,dotsep=0.7pt](0,0)(0,4.061)
\psline[linewidth=1pt](-.007,1.01)(0,1.015)
\psdots[dotstyle=*,dotsize=2pt](0,1.015)
\psline[linewidth=1pt,linestyle=dashed](-.007,3.04)(0,3.046)
\psdots[dotstyle=*,dotsize=2pt](0,3.046)
\uput[dr](.495,.751){$\scriptstyle \alpha^{}_{\times}$}
\uput[r](0,-1.625){$\scriptstyle \gamma^+_{\times} \sim\, \gamma^+_0$}
\uput[u](2.166,-.268){$\scriptstyle \gamma^-_{\times} \sim\, \gamma^-_0$}
\uput[d](-1.083,.134){$\scriptstyle q^-_{\delta , n}$}
\uput[ur](0,1.015){$\scriptstyle q^+_{\delta , n}$}
\uput[d](-3.25,.402){$\scriptstyle \Phi^{}_{\delta} (q^-_{0 , n})$}
\uput[ur](0,3.046){$\scriptstyle \Phi^{}_{\delta} (q^+_{0 , n})$}
\uput[ur](.928,1.408){$\scriptstyle \L^{\cap}_{\W_{\times}}$}
\uput[ur](-3.404,1.944){$\scriptstyle \L^-_{\W_{\times}}$}
\uput[dr](.928,5.469){$\scriptstyle \L^+_{\W_{\times}}$}
\end{pspicture}
\caption{The extension of a segment $\check{c}^-_{\delta}\!\!\mid_{[0,n \check{l}_{\delta}]} \dtimes \id^{}_{\EE}$ and the quasi-isometric image within the wall $\W^{}_{\times , n} \subset X_{\delta}$.}
\label{fig:Wx}
\end{figure}
We now choose a fringe $\L^+_{\times , n}$ covering $E^{}_{\times} \cap C^+_{\times}$ of distance less than $\smash{\sqrt{2}}$ above $\smash{q^-_{\delta , n}}$ with respect to $\gamma^-_{\times}$.
The point $\smash{\Phi^{}_{\delta} (q^-_{0 , n})}$ then lies at a distance $\smash{n ( \check{l}^{}_{\delta} - \check{l}^{}_0 )}$ below $q^-_{\delta , n}$.
By choosing points $\smash{x^{}_{ 0 , n } \in \pi^{-1} (p^+_0) \subset X^{}_0}$ and $\smash{x^{}_{ \delta , n } = \Phi^{}_{\delta} (x^{}_{ 0 , n }) \in X^{}_{\delta}}$ we construct a corresponding sequence of geodesic segments in $\smash{\check{X}^+_0}$ and $\smash{\check{X}^+_{\delta}}$ such that their endpoints $\smash{q^+_{\delta , n}}$ and $\smash{\Phi^{}_{\delta} (q^+_{0 , n})}$ come to lie on the same axis of $\smash{\gamma^+_{\times}}$ in the fringe $\smash{\L^+_{\times , n}}$ defined above, having the same distances to $\smash{\L^{\cap}_{\times , n} = \L^-_{\times , n} \cap \L^+_{\times , n}}$ as $\smash{q^-_{\delta , n}}$ and $\smash{\Phi^{}_{\delta} (q^-_{0 , n})}$.
Since the fringes $\L^-_{\times , n}$ and $\L^+_{\times , n}$ in $\W^{}_{\times , n}$ are orthogonal (see Figure \ref{fig:Wx}),
\begin{align}
c^{}_{ 0 , n } & = \overline{x^{}_0 q^-_{ 0 , n}} \cup \overline{q^-_{ 0 , n} q^+_{ 0 , n}} \cup \overline{q^+_{ 0 , n} x^{}_{ 0 , n}} \\
\shortintertext{and}
c^{}_{ \delta , n } & = \overline{\Phi^{}_{\delta} (x^{}_0) q^-_{\delta , n}} \cup \overline{q^-_{\delta , n} q^+_{\delta , n}} \cup \overline{q^+_{\delta , n} \Phi^{}_{\delta} (x^{}_{ 0 , n})} \\
	& = \overline{x^{}_{\delta} q^-_{\delta , n}} \cup \overline{q^-_{\delta , n} q^+_{\delta , n}} \cup \overline{q^+_{\delta , n} x^{}_{ \delta , n}}
\end{align}
are geodesic segments.

\subsection{Linear divergence}
\label{sec:cex:non_stability}
Concluding this section we will prove that the quasi-isometry $\Phi_{\delta}$ violates any sublinear estimate as in \cite[Lemma 2.5]{crkl_geod} and does not extend continuously to the virtual boundary of $X_0$.

By our choice of the fringe $\smash{\L^{\cap}_{\times , n}}$ we have
\begin{align}
l (c^{}_{\delta , n}) & = 2 \sqrt{2} \; l \bigl( \overline{\check{c}^-_{\delta} (0) \check{q}^-_{\delta , n}} \dtimes \id^{}_{\EE} \bigr) + d_{X^{}_{\delta}} \bigl( q^-_{\delta,n} , q^+_{\delta,n} \bigr) \\
	& \leqslant 2 \sqrt{2} ( n \check{l}^{}_{\delta} + 1 ) + 2 \; ,\\
\shortintertext{and}
d^{}_{X^{}_0} \bigl( q^-_{0,n} , q^+_{0,n} \bigr) & = d^{}_{X^{}_{\delta}} \bigl( \Phi^{}_{\delta} ( q^-_{0,n} ) , \Phi^{}_{\delta} ( q^+_{0,n} ) \bigr) \\
	& = \sqrt{2} \, d^{}_{X^{}_{\delta}} \bigl( \Phi^{}_{\delta} ( q^-_{0,n} ) , \L^{\cap}_{\times , n} \bigr) \\
	& = \sqrt{2} \, \bigl( d^{}_{X^{}_{\delta}} \bigl( \Phi^{}_{\delta} ( q^-_{0,n} ) , q^-_{\delta , n} \bigr) + d^{}_{X^{}_{\delta}} \bigl( q^-_{\delta , n} , \L^{\cap}_{\times , n} \bigr) \bigr) \\
	& \leqslant \sqrt{2} n ( \check{l}^{}_{\delta} - \check{l}^{}_0 ) + 2 \\
\shortintertext{yields}
l (c^{}_{0, n}) & = 2 \sqrt{2} \; l \bigl( \overline{\check{c}^-_0 (0) \check{q}^-_{0 , n}} \dtimes \id^{}_{\EE} \bigr) + d^{}_{X^{}_0} \bigl( q^-_{0,n} , q^+_{0,n} \bigr) \\
	& \leqslant 2 \sqrt{2} ( n \check{l}^{}_0 + 1 ) + \sqrt{2} n ( \check{l}^{}_{\delta} - \check{l}^{}_0 ) + 2 \\
	& = n \sqrt{2} ( \check{l}^{}_{\delta} + \check{l}^{}_0 ) + 2 ( 1 + \sqrt{2} ) \; .
\end{align}
The midpoint $m^{}_{\varepsilon , n}$ of $c^{}_{ \varepsilon , n }$ is also midpoint of $\overline{q^-_{\varepsilon , n} q^+_{\varepsilon , n}}\,$; hence $\Phi^{}_{\delta} (m^{}_{0 , n})$ is the midpoint of $\overline{\Phi^{}_{\delta} (q^-_{0 , n}) \Phi^{}_{\delta} (q^+_{0 , n})}$.
Note that $\proj^{}_{c^{}_{\delta , n}} \bigl(\Phi^{}_{\delta} (m^{}_{0 , n}) \bigr) = m^{}_{\delta , n}\,$.
Hence
\begin{equation}
d_{X^{}_{\delta}} \bigl( \Phi^{}_{\delta} ( m^{}_{0,n} ) , c^{}_{ \delta , n } \bigr)
	= \frac{1}{\sqrt{2}} d_{X^{}_{\delta}} \bigl( \Phi^{}_{\delta} ( q^-_{0,n} ) , q^-_{ \delta , n } \bigr) \\
	= \frac{n}{\sqrt{2}} ( \check{l}^{}_{\delta} - \check{l}^{}_0 ) \; .
\end{equation}
Thereby the distance between the paths $\smash{\Phi_{\delta} \circ c_{0,n}}$ and $\smash{c_{\delta,n}}$ within the wall $\smash{\W_{\times , n}}$ increases linearly in $n$. Furthermore,
\begin{equation}
\frac{d_{X^{}_{\delta}} \bigl( \Phi^{}_{\delta} ( m^{}_{0,n} ) , c^{}_{ \delta , n } \bigr)}{1 + \frac{1}{2} l (c^{}_{ 0 , n }) }
	\geqslant \frac{ \frac{n}{\sqrt{2}} ( \check{l}^{}_{\delta} - \check{l}^{}_0 ) }{ 1 + \frac{n}{\sqrt{2}} ( \check{l}^{}_{\delta} + \check{l}^{}_0 ) + 1 + \sqrt{2} } \\
	\geqslant \frac{\check{l}^{}_{\delta} - \check{l}^{}_0}{\check{l}^{}_{\delta} + \check{l}^{}_0 + 5}
\end{equation}
implies that there is no real function $\theta$ vanishing at infinity such that for all $x,y \in X_0$ and $z \in \overline{xy}$ the sublinear estimate
\begin{equation} \label{eq:sublinear_estimate}
d^{}_{X_{\delta}} \bigl( \Phi^{}_{\delta} (z) , \overline{\Phi^{}_{\delta} (x) \Phi^{}_{\delta} (y)} \bigr) \leqslant \bigr( 1 + d^{}_{X_0}(z,x) \bigr) \; \theta ( d^{}_{X_0}(z,x) )
\end{equation}
holds.
Consequently, the generalized geometric data carry insufficient geometric information as to prove a generalization of \cite[Theorem 1.3]{crkl_geod}.
Moreover, geodesic rays or segments in $X^{}_{\delta}$ cannot be modeled on templates as in \cite[Theorem 5.1]{crkl_geod}.

\begin{prop} \label{prop:non_continuity}
The quasi-isometry $\Phi_{\delta}$ does not extend continuously to the virtual boundary of $X_0$.
\end{prop}

\begin{pr}
The sequence $x^{}_{\varepsilon , n}$ converges to a point $\vb{x}^{}_{\varepsilon}$ in $\partial_{\infty} X^{}_{\varepsilon}$ since for any $k$ and all $n \geqslant k$ the geodesic segments $c^{}_{\varepsilon , k}$ and $c^{}_{\varepsilon , n}$, reparametrised to arclength, coincide up to time $\smash{\sqrt{2} k \check{l}^{}_{\varepsilon}}$.
Hence for $\smash{ x'^{}_{\varepsilon , k} = \overline{x^{}_{\varepsilon , 0} \vb{x}^{}_{\varepsilon}} (\sqrt{2} k \check{l}^{}_{\varepsilon})}$ we have $\smash{x'^{}_{\varepsilon , k} = c^{}_{\varepsilon , n} (\sqrt{2} k \check{l}^{}_{\varepsilon})}$ for any $n \geqslant k$.
Regarding $( \check{c}^-_{\delta} (0) , 0 )$ as the point of origin in $\smash{\minset_{X_{\delta}} (\gamma^-_{\times}) \cong \check{X}^-_{\delta} \times \EE}$ let $\smash{y^{}_{\delta, n}}$ denote the point $\Phi^{}_{\delta} (x'^{}_{0 , n}) = ( \check{o}^-_{4k} , k \check{l}^{}_0 ) \in \bigl( \check{\L}^-_{4k} \cap \check{\L}^+_{4k} \bigr) \times \EE$ with respect to the splitting.
For the constant speed parametrization of the segment $\smash{\overline{x^{}_{\delta} y^{}_{\delta , n}}}$ on the interval $\smash{[0 , n\check{l}^{}_{\delta}]}$ we have $\overline{x^{}_{\delta} y^{}_{\delta , n}} (k \check{l}^{}_{\delta}) = y^{}_{\delta , k}$ for all $k \leqslant n$.
Therefore the sequence $\smash{y^{}_{\delta , n}}$ converges to a point $\smash{\vb{y}^{}_{\delta} \in \partial^{}_{\infty} X^{}_{\delta}}$.
Since $\smash{c^-_{\delta , n}}$ and $\smash{\Phi^{}_{\delta} \circ c^-_{0 , n}}$ are parallel within $\smash{\S^-_{\times , n} \cup \W^{}_{\times , n} \cup \S^+_{\times , n}}$, we have
$\smash{d_{X^{}_{\delta}} \bigl( x'^{}_{\delta , n} , y^{}_{\delta , n} \bigr) = d_{X^{}_{\delta}} \bigl( q^-_{\delta , n} , \Phi^{}_{\delta} (q^-_{0 , n}) \bigr) = n (\check{l}^{}_{\delta} - \check{l}^{}_0 ) \geqslant (\check{l}^{}_{\delta} - \check{l}^{}_0 ) > 0}$.
Hence $\vb{x}_{\delta} \neq \vb{y}_{\delta}$, whereas any continuous extension of $\Phi_{\delta}$ in $\vb{x}_0$ would yield equality.
\qed
\end{pr}


\section{Retrieval of geometric data}
\label{sec:salvage}
After having demonstrated that equivalence of the geometric data of two admissible rank $3$ actions is not sufficient to ensure boundary stability we will discuss as to what extent the homeomorphism class of the virtual boundary determines the equivalence class of the geometric data.
It will be shown that the geometric data of two admissible rank $3$ actions are, for the most part, equivalent if an equivariant quasi-isometry extends continuously to the virtual boundaries.

In what follows, let $\Gamma \act X$ be an admissible action of rank 3; i.e., $\Gamma$ acts geometrically on a locally compact Hadamard space $X$ and there is a $\Gamma$-tree $\T$ with at least one edge such that the following conditions are satisfied:
\begin{enumerate}
\item Any vertex stabilizer $\Gamma_v$ is an extension of a non-elementary hyperbolic group by a free abelian group  $\Lambda_v < \Gamma_v$ of rank $2$.
	\label{rem_2:adm_tr_gr:vertex_groups}
	
\item Any edge stabilizer $\Gamma_e$ is a finite extension of a free abelian group $\Lambda_e$ of rank $3$ and the isotropy groups of two edges incident to the same vertex are incommensurable.
	\label{rem_2:adm_tr_gr:edge_groups}

\item For any edge $e$ the subgroup generated by $\Lambda_{\partial^- e} \Lambda_{\partial^+ e}$ has finite index\footnote{This is equivalent to $\comm_{\Gamma} (\Lambda_v) = \Gamma_v$ for all $v \in \V$.} in $\Lambda_e$.
	\label{rem_2:adm_tr_gr:vertex_union}

\item The families $\{ \Lambda_v \}_{v \in \V\!\!\T}$ and $\{ \Lambda_e \}_{e \in \E\!\T}$ are $\Gamma$-equivariant.
	\label{rem_2:adm_tr_gr:equivariance}
\end{enumerate}

\subsection{Generalized blocks}
\label{sec:salvage:gen_sfp}
First of all, we construct a family of subspaces in $X$, each splitting off a factor that allows for an admissible action of rank 2, to which we can apply Theorem 1.3 in \cite{crkl_geod}.

\begin{prop} \label{prop:reduced_adm_action}
For any edge $e^*$ of $\T$ there exists an infinite cyclic subgroup $\Lambda$ in $\smash{\Lambda_{\del^- e^*} \cap \Lambda_{\del^+ e^*}}$ such that the action of $\check{\Gamma} = \smash{\rquotient{\normalizer_{\Gamma} (\Lambda)}{\Lambda}}$ on the non-Euclidean%
\footnote{$\check{X}$ does indeed contain the non-Euclidean convex subspace $\overline{Y}_{\partial^- e^*}$.}
factor of $\smash{\minset_X (\Lambda) \cong \check{X} \times \EE}$ is admissible of rank $2$ and its Bass-Serre tree $\smash{\check{\T}}$ is a subtree of $\T$.
\end{prop}

\begin{pr}
For each edge $e$ of $\T$ the subgroups $\smash{\Lambda_{\del^- e}}$ and $\smash{\Lambda_{\del^+ e}}$ in $\Lambda_e$ are incommensurable as per condition \ref{rem_2:adm_tr_gr:vertex_union}; hence their intersection $\smash{\Lambda^{\cap}_e = \Lambda_{\del^- e} \cap \Lambda_{\del^+ e}}$ is infinite cyclic.

The subgraph $\smash{\check{\T}}$ defined by $\{ e \in \E\!\T \mid \Lambda^{\cap}_{e} \sim \Lambda^{\cap}_{e^*} \}$ is a tree containing at least the one edge $e^*$ and upon which $\Gamma^* = \comm_{\Gamma} (\Lambda^{\cap}_{e^*})$ acts with finite quotient $\smash{\check{\G}  = \rquotient{\check{\T}\!\!}{\Gamma^*}}$:
For an edge $e$ in $\check{\T}$ let $(v_0 e_1 v_1 \ldots v_{n-1} e_n v_n )$ be the minimal reduced path from $e = e_1$ to $e^{*} = e_n$ in $\T$.
For any $1 < j < n$ we then have
\begin{equation}
\Lambda^{\cap}_{e^{*}} \sim \, \Lambda^{\cap}_e \cap \, \Lambda^{\cap}_{e^{*}} < \Gamma^{}_{v^{}_0}\!\cap \, \Gamma^{}_{v^{}_n}\!< \Gamma^{}_{e^{}_{j-1}}\!\!\cap \, \Gamma^{}_{e^{}_j}\!\cap \, \Gamma^{}_{e^{}_{j+1}}\!\!\sim \, \Lambda^{}_{v^{}_{j-1}}\!\cap \, \Lambda^{}_{v^{}_j} = \, \Lambda^{\cap}_{e^{}_j} \; .
\end{equation}
Hence $\smash{\check{\T}}$ is connected.
Moreover, $\smash{\Lambda^{\cap}_{\gamma.e} = {}^\gamma \Lambda^{\cap}_{e} \sim {}^\gamma \Lambda^{\cap}_{e^*} \sim \, \Lambda^{\cap}_{e^*}}$ holds for any $e \in \E\!\check{\T}$ and $\gamma \in \Gamma^*$.
Thus $\smash{\check{\T}}$ is $\Gamma^*$-invariant.
Let $e$ be an edge of $\check{\T}$ and, on the other hand, $\gamma \in \Gamma$ such that $\gamma . e \in \E\!\check{\T}$
then $\Lambda^{\cap}_{e^*} \sim \, \Lambda^{\cap}_{\gamma.e} = {}^\gamma \Lambda^{\cap}_{e} \sim {}^\gamma \Lambda^{\cap}_{e^*}$ and $\gamma$ is actually contained in $\Gamma^*$.
Therefore the canonical map $\rquotient{\check{\T}\!\!}{\Gamma^*} \rightarrow \rquotient{\T\!\!}{\Gamma}$ is injective; hence the quotient $\rquotient{\check{\T}\!\!}{\Gamma^*}$ is finite.

Let $\gr{S}$ be a connected fundamental domain%
\footnote{This is to be understood in the strict sense, so $\gr{S}$ might not be a subgraph of $\T$.}
of the $\Gamma^*$-action on $\check{\T}$.
Then $\Lambda = \bigcap_{e \in \E\!\gr{S}} \Lambda^{\cap}_e$ is a cyclic normal subgroup of $\Gamma^*$:
Since $\gr{S}$ is finite, $\Lambda$ and $\smash{\Lambda^{\cap}_{e^*}}$ are commensurable and $\Lambda$ is infinite cyclic.
After choosing a maximal subtree $\gr{S}^{}_0$ in $\gr{S}$ the set $\Sigma = \bigl\{ \tau^{}_e \mid e \in \E\!\gr{S} \setminus \E\!\gr{S}^{}_0 \bigr\} \cup \bigcup_{v \in \V\!\!\gr{S}}\nolimits \Gamma^*_v$ forms a system of generators%
\footnote{For an $e \in \E\!\gr{S} \setminus \E\!\gr{S}^{}_0$ the element $\tau^{}_e$ acts as a translation of a geodesic in $\T$ through $e$.}
of $\Gamma^*$; see section \ref{sec:fund_grp}.
It therefore suffices to show that $\Sigma \subset \normalizer^{}_\Gamma (\Lambda)$:
Let $\xi$ be a generator of $\Lambda$ and $\sigma \in \Sigma$.
Since $\Lambda$ and $\smash{{}^\sigma \Lambda}$ are commensurable and the translation distance of $\xi$ is positive and invariant under conjugation, there exists a $k \in \NN$ such that $\smash{{}^\sigma \xi^k = \xi^{\pm k}}$.
In case $\sigma \in \Gamma^*_v$ for some vertex $v$ in $\gr{S}$, we have ${}^\sigma \xi \in {}^\sigma \Lambda^{}_v = \Lambda^{}_{\sigma . v} = \Lambda^{}_v$.
If otherwise $\sigma = \tau^{}_e$ for an $\smash{e \in \E\!\gr{S} \setminus \E\!\gr{S}^{}_0}$, the vertex $v = \partial^+ ( \tau^{}_e . e )$ lies in $\gr{S}$, and we have $\smash{{}^{\tau^{}_e} \xi \in {}^{\tau^{}_e} \Lambda^{\cap}_e = \Lambda^{\cap}_{\tau^{}_e . e} < \Lambda^{}_v}$.
In either case, there exists a vertex $v \in \V\!\gr{S}$ such that $\{ \xi , {}^\sigma \xi \} \subset \Lambda^{}_v$.
Therefore $\xi$ and $\eta = ( {}^\sigma\xi ) \xi^{-1}$ commute.
In particular, $\eta^k \xi^k = ( \eta \xi )^k = {}^\sigma\xi^{k} = \xi^{\pm k}$; hence ${}^\sigma \xi = \xi^{\pm 1}$ and ${}^\sigma\Lambda = \Lambda$.%

It remains to show that the induced action of $\check{\Gamma} = \smash{\rquotient{\normalizer^{}_{\Gamma} (\Lambda)}{\Lambda}}$ on the non-Euclidean factor of $\smash{\minset_X (\Lambda) \cong \check{X} \times \EE}$ is admissible of rank 2.
First of all, $\check{\Gamma}$ acts geometrically on $\check{X}$.
From the discussion above it follows in particular that $\Lambda$ is contained in the isotropy subgroup of any simplex of $\check{\T}$.
Hence $\check{\Gamma}$ acts on $\check{\T}$ with $\smash{\check{\Gamma}_a = \rquotient{\bigl( \Gamma_a \cap \normalizer_{\Gamma} (\Lambda) \bigr)}{\Lambda}}$ for any edge or vertex $a$ of $\check{\T}$.
Let $v$ be a vertex of $\check{\T}$.
With $\Lambda < \Lambda_v$ we deduce from $\Lambda_v < \centralizer_{\Gamma} (\Lambda_v) < \centralizer_{\Gamma} (\Lambda) \sim \normalizer^{}_{\Gamma} (\Lambda) = \Gamma^*$ that, for one, $\Lambda^{}_v$ is a normal subgroup of $\Gamma^*_v$ and, secondly, that $\Gamma^*_v$ has finite index in $\Gamma^{}_v$.
Since $\rquotient{\Lambda_v}{\Lambda}$ is virtually infinite cyclic we find therein an infinite cyclic characteristic subgroup $\check{\Lambda}_v$ of finite index.
The induced action of $\check{\Delta}^{}_v = \rquotient{\check{\Gamma}^{}_v}{\check{\Lambda}^{}_v}$ on the hyperbolic space $\overline{Y}^{}_v$ is geometric because $\check{\Delta}^{}_v$ is a finite extension of $\rquotient{\Gamma^*_v}{\Lambda^{}_v}$.
Since $\overline{Y}^{}_v$ has more than two ends, $\check{\Delta}^{}_v$ is non-elementary hyperbolic and \ref{rem_2:adm_tr_gr:vertex_groups} holds.
A similar argument shows $\Lambda^{}_e < \Gamma^*_e$.
Consequently $\check{\Gamma}^{}_e$ contains the virtually free abelian rank 2 subgroup $\rquotient{\Lambda^{}_e}{\Lambda}$ with finite index; hence we have shown \ref{rem_2:adm_tr_gr:edge_groups}.
Both commensurability conditions follow readily from the fact that two subgroups in $\check{\Gamma}$ are commensurable if and only if this holds for their lifts.
The family is obviously $\check{\Gamma}$-equivariant.
\qed
\end{pr}

We hereby obtain, for each edge $e$ of $\T$, a locally compact Hadamard space $\check{X} (e)$ with an admissible rank $2$ action by a group $\check{\Gamma} (e)$ with Bass-Serre tree $\check{\T} (e)$.

\subsection{Stability under quasi-isometries}
\label{sec:salvage:qi_stab}
Let $\Gamma \act X$ and $\Gamma \act X'$ be admissible actions of rank $3$ and $\Phi$ an equivariant quasi-isometry from $X$ to $X'$ which extends to an equivariant homeomorphism $\virtbnd \Phi$ of the virtual boundaries.

For any edge $e$ of $\T$ the actions $\check{\Gamma} (e) \act \check{X} (e)$ and $\check{\Gamma} (e) \act \check{X}' (e)$ are, according to Proposition \ref{prop:reduced_adm_action}, cocompact.
Consequently $\check{X} = \check{X} (e)$ and $\check{X}' = \check{X}' (e)$ are equivariantly quasi-isometric as well.
We are about to deduce the equivalence of the geometric data of the actions $\Gamma \act X$ und $\Gamma \act X'$ by applying \cite[Theorem 1.3]{crkl_geod} to the induced rank $2$ actions on the subspaces $\check{X}$ and $\check{X}'$ respectively.
To that end, we have to show that any equivariant quasi-isometry between these subspaces can be extended to a homeomorphism of their virtual boundaries.

Let $\Lambda = \Lambda (e)$ denote the infinite cyclic subgroup of $\Lambda^{\cap}_e$ constructed in Proposition \ref{prop:reduced_adm_action}.
It follows from Corollary \ref{cor:min_mr_fin_hd} that $\Phi \bigl( \minset^{}_X (\Lambda) \bigr)$ and $\minset^{}_{X'} (\Lambda)$ have finite Hausdorff distance. 
The image of $\virtbnd\minset^{}_X (\Lambda)$ under $\virtbnd\Phi$ is therefore contained in $\virtbnd\minset^{}_{X'} (\Lambda)$.
The product decomposition $\minset^{}_X (\Lambda) \cong \check{X} \times \EE$ yields a homeomorphism $\virtbnd\minset^{}_X (\Lambda) \cong \Sigma\, \virtbnd\check{X}$; an analogous statement holds in $X'$.
Since the real fibres in $\minset_X (\Lambda)$ and $\minset_{X'} (\Lambda)$ are the convex closures of the $\Lambda$-orbits, the equivariance of $\Phi$ implies that the image of any $\EE$-fibre under $\Phi$ lies in uniformly bounded Hausdorff distance to an $\EE$-fibre.
With respect to the aforementioned homeomorphism the poles of the suspension $\Sigma\,\smash{\virtbnd\check{X}}'$ do therefore correspond to the images of the poles of $\Sigma\,\virtbnd\check{X}$ under $\virtbnd\Phi$.

In what follows, we will denote the coordinates of an $x \in \minset^{}_X (\Lambda)$ with respect to the splitting $\minset^{}_X (\Lambda) \cong \check{X} \times \EE$ by $(\check{x} , t^{}_x)$.
After fixing basepoints $x^{}_0 \in \minset^{}_X (\Lambda)$ and $y^{}_0 = \Phi (x^{}_0)$ we may assume that $x^{}_0 = (\check{x}^{}_0 , 0 )$ and $y^{}_0 = (\check{y}^{}_0 , 0 )$.

\begin{lem} \label{lem:gen_sfp_qi_reduction}
The map $\smash{\check{\Phi} = \proj_{\check{X}'} \circ \proj_{\minset_{X'} (\Lambda)} \circ \, \Phi|_{\check{X}}}$ is a $\check{\Gamma}$-equivariant quasi-isometry and induces a $\check{\Gamma}$-equivariant homeomorphism from $\virtbnd\check{X}$ onto $\virtbnd\check{X}'$.
\end{lem}

\begin{pr}
Whenever $\Phi$ satisfies an $(L,A)$-Lipschitz condition, the same holds for $\check{\Phi}$ because any projection has Lipschitz constant $1$.
Since $\Lambda$ acts cocompactly on the $\EE$-fibres, there exists a constant $C$ such that for any $x \in \check{X} \times \{ 0 \}$ there is a $\gamma^{}_x \in \Lambda$ with $d^{}_{X'} \left( \check{\Phi} (\check{x}) , \gamma^{}_x . \Phi (x) \right) \leqslant C$.
For all $x , y \in \check{X} \times \{ 0 \}$ we therefore have $\smash{d_{\check{X}'} \left( \check{\Phi} (\check{x}) , \check{\Phi} (\check{y}) \right) \geqslant L^{-1} \smash{d_{\check{X}}} \left( \check{x} , \check{y} \right) - A - 2C}$.
%
%
It is easy to check that $\smash{\check{\Phi} (\check{X})}$ is quasi-dense in $\smash{\check{X}'}$. The $\smash{\check{\Gamma}}$-equivariance follows immediately from the $\Gamma$-equivariance of $\Phi$.

Since $\Phi$ extends continuously to $\virtbnd\check{X}$ by assumption, it suffices to show that for a sequence $y_n = ( \check{y}_n , t_n )$ in $\smash{\minset_{X'} (\Lambda)}$ with $\lim_{n \rightarrow \infty} y_n = \vb{y} \in \virtbnd\Phi \bigl( \virtbnd\check{X} \bigr)$ the sequence $\check{y}_n$ in $\smash{\check{X}}'$ converges to an $\check{\vb{y}} \in \virtbnd\smash{\check{X}}'$.
For any such sequence and for any $R>0$ the sequence $n \mapsto \smash{y^R_n = \proj_{\cball{R}{y_0}} (y_n)}$ is convergent as well.
Since $\vb{y}$ is not a pole of the suspension and, for suitable $n(R)$, the sequence $y^R_{n(R)}$ does also converge  to $\vb{y}$ as $R$ tends to infinity, we can, for any $r>0$, choose $R>0$ and $N$ such that $\smash{d_{\check{X}'} \bigl( \check{y}^{}_0 , \proj_{\check{X}'} (y^R_{N+k}) \bigr) > r}$ holds for all $k$.
With $d_{X'} ( y_0 , y_n ) = \sqrt{\smash{d_{\check{X}'}} ( \check{y}_0 , \check{y}_n )^2 + \smash{\norm{t_n}'}^2 }$ tending to infinity, $d_{\check{X}'} ( \check{y}^{}_0 , \check{y}^{}_n )$ is unbounded as well.
Thereby we have for sufficiently large $N$ that $\smash{\proj_{\cball{r}{\check{y}^{}_0}} (\check{y}^{}_{N+k}) = \proj_{\cball{r}{\check{y}^{}_0}} \bigl( \proj_{\check{X}'} (y^R_{N+k}) \bigr)}$ for all $k$.
Since projections onto convex subspaces are distance non-increasing, $\smash{\check{y}^r_n = \proj_{\cball{r}{\check{y}^{}_0}} (\check{y}^{}_n)}$ is a Cauchy sequence.
Hence $\check{y}^{}_n$ converges in $\check{X}$ to a $\check{\vb{y}} \in \virtbnd\check{X}'$. 

Let $c$ denote the parametrization of $\smash{\overline{y_0 \vb{y}}}$ by arclength.
Then $\check{c} = \smash{\proj_{\check{X}'}} \circ \, c$ is a parametrization of the ray $\smash{\overline{\check{y}_0\check{\vb{y}}}}$ with constant speed $\sqrt{1 - \smash{d_{X'}} (\check{X}' , c (1))^2}$.
Now, let $\check{c}$ denote the parametrization of $\overline{\check{y}_0\check{\vb{y}}} \subset \check{X}'$ by arclength and be $\vb{x}_1 , \vb{x}_2 \in \virtbnd\check{X}$ with $\virtbnd\check{\Phi} (\vb{x}_2) = \virtbnd\check{\Phi} (\vb{x}_1) = \check{\vb{y}}$.
Then there are $\theta_1 , \theta_2 \in \bigl( -\frac{\pi}{2} , \frac{\pi}{2} \bigr)$ such that for $i \in \{ 1,2 \}$ the geodesic ray $c_i : t \mapsto ( \check{c} ( t \cos \theta_i ) , t \sin \theta_i )$ is a representative of $\virtbnd\Phi (\vb{x}_i)$.
After choosing an equivariant quasi-inverse of $\Phi$ we find sequences $x^{}_n$ in $\minset^{}_X (\Lambda)$ and $\gamma^{}_n$ in $\Lambda$ such that $\lim^{}_{n \rightarrow \infty} x^{}_n = \vb{x}^{}_1$ and $\lim^{}_{n \rightarrow \infty} \gamma^{}_n . x^{}_n = \vb{x}^{}_2$.
Since $\vb{x}^{}_1 , \vb{x}^{}_2 \in \virtbnd\check{X}$ and $\proj^{}_{\check{X}} (x^{}_n) = \proj^{}_{\check{X}} (\gamma^{}_n . x^{}_n)$ holds for all $n$, we have $\vb{x}^{}_1 = \vb{x}^{}_2$.
Hence $\virtbnd\check{\Phi}$ is injective.

A similar argument shows that $\virtbnd\check{\Phi}$ is surjective and its inverse is continuous.
\qed
\end{pr}

\subsection{Equivalence of geometric data}
\label{sec:salvage:geom_data_ret}
Let $v$ be a vertex of $\T$ and $e$ any edge incident to it.
Then, by \cite[Theorem 1.3]{crkl_geod} and Lemma \ref{lem:gen_sfp_qi_reduction}, the geometric data in the sense of \textsc{Croke} and \textsc{Kleiner} of the actions $\check{\Gamma} (e) \act \check{X} (e)$ and $\check{\Gamma} (e) \act \check{X}' (e)$ are equivalent, i.e., for any vertex $w$ in $\check{\T}$ there exist positive constants $\lambda(w)$ and $\mu(w)$ such that $\smash{\check{\mls}}'_w = \lambda(w) \smash{\check{\mls}}_w$ and $\smash{\check{\tau}}'_w = \mu(w) \check{\tau}_w$.

Since $\Lambda = \Lambda (e)$ is central in $\Lambda^{}_v$, we may regard the non-Euclidean factor $\overline{Y}^{}_v$ of $\smash{\minset_{X} (\Lambda_v)}$ as a subspace in $\check{X}$.
The construction of the action $\check{\Gamma} \act \check{X}$ shows that there is a subgroup $\tilde{\Lambda}^{}_v$ of finite index in $\Lambda^{}_v$ such that $\rquotient{\tilde{\Lambda}^{}_v}{\Lambda} = \check{\Lambda}^{}_v$.
We may therefore assume that $\smash{\overline{Y}^{}_v}$ is contained in the non-Euclidean factor $\smash{\check{Y}^{}_v}$ of $\smash{\minset_{\check{X}} (\check{\Lambda}^{}_v)}$.
Now, fix a $\gamma \in \Gamma^{}_v$.
Since $\centralizer^{}_{\Gamma} (\Lambda^{}_v)$ is a subgroup of $\normalizer^{}_{\Gamma} (\Lambda)$ and has finite index in $\Gamma^{}_v$, $\Gamma_v \cap \normalizer_{\Gamma} (\Lambda)$ has finite index in $\Gamma^{}_v$ as well.
Hence there is a $k \in \NN$ for which $\check{\gamma}^k = \gamma^k \Lambda \in \check{\Gamma}^{}_v$, so $\mls^{}_v (\gamma) = k^{-1} \check{\mls}^{}_v (\check{\gamma}^k)$ holds.
The map $\mls_v$ defined on $\Gamma_v$ is therefore uniquely determined by $\check{\mls}_v$.
Consequently the geometric data $\mls$ and $\mls'$ are equivalent. 

Furthermore, the geometric data $\Tau_v$ and $\Tau'_v$ coincide for all those vertices $v$ which are contained in the intersection of distinct subtrees $\check{\T} (e_1)$ and $\check{\T} (e_2)$:
Let $v$ be a vertex of $\T$ with adjacent edges $e^{}_1 , e^{}_2$ such that $\smash{\Lambda^{\cap}_{e_1}}$ and $\smash{\Lambda^{\cap}_{e^{}_2}}$ are not commensurable%
\footnote{In case there were no such vertex, there existed an infinite cyclic subgroup $\Lambda$ in $\Gamma$ with finite index centralizer in $\Gamma$ and the Hausdorff distance between $X$ and $\minset^{}_X (\Lambda) \cong \check{X} \times \EE$ were finite.}.
By applying Proposition \ref{prop:reduced_adm_action} to the edges $e^{}_1$ and $e^{}_2$ we obtain for each $i$ an infinite cyclic subgroup $\smash{\Lambda_i < \Lambda^{\cap}_{e_i}}$ such that the induced action of the quotient $\smash{\check{\Gamma}_i = \rquotient{\normalizer_{\Gamma} (\Lambda_i)}{\Lambda_i}}$ on the non-Euclidean factor $\check{X}_i$ of $\minset_X (\Lambda_i)$ is admissible.
Let $\check{\Lambda}_{v,i}$ denote the infinite cyclic subgroup of $\rquotient{\Lambda_v}{\Lambda_i}$ we constructed in the proof of Proposition \ref{prop:reduced_adm_action}.
For each $i$ we choose a generator of $\check{\Lambda}_{v,i}$.
Since $\centralizer_{\Gamma} (\Lambda_v)$ is a subgroup in each $\normalizer_{\Gamma} (\Lambda_i) \cap \Gamma_v$, the second geometric data in the sense of \cite{crkl_geod} then correspond to $\Lambda_i$-invariant maps $\tau_i$ from $\centralizer_{\Gamma} (\Lambda_v)$ to $\RR$.

As in section \ref{sec:geom_data_def} we regard the two-dimensional Euclidean factor $\EE_v$ of $\smash{\minset_X (\Lambda_v)}$ as the real vector space $\smash{\Lambda_v \otimes_{\ZZ} \RR}$ with the metric $<.\,,.>$ induced by the quadratic form $\zeta \mapsto \abs{\zeta}^2_X$.
A pair of generators $\xi_1$ and $\xi_2$ of $\Lambda_1$ and $\Lambda_2$ resp.\ is then linearly independent because $\Lambda_1$ and $\Lambda_2$ are not commensurable.
For each $i$ let $w_i \in \EE_v$ denote a unit vector in the orthogonal complement of $\xi_i$ such that for a lift $\zeta_i$ of the generator of $\check{\Lambda}_{v,i}$, chosen to define $\tau_i$, the bases $(w_i , \xi_i)$ and $(\zeta_i , \xi_i)$ have the same orientation.
Thereby we may identify $\tau_i$ with the linear functional $\gamma \mapsto < \gamma_{\EE} , w_i >$ on $\EE_v$, where $\gamma_{\EE}$ denotes the translation vector of $\gamma$ on $\EE_v$.
Especially $\tau_i (\xi_i) = 0$ and $\tau_i (\xi_j) \neq 0$ for $j \neq i$ since $\xi_1$ and $\xi_2$ are linearly independent.
The preceding discussion holds analogously for $X'$ and in what follows, primed labels are to be understood referring to it.

The virtual boundaries of $\check{X}_1$ and $\smash{\check{X}}'_1$ as well as $\check{X}_2$ and $\smash{\check{X}}'_2$ are equivariantly homeomorphic by Lemma \ref{lem:gen_sfp_qi_reduction}.
Thus, according to \cite[Theorem 1.3]{crkl_geod}, there exists for each $i$ a positive constant $\mu_i$ such that $\tau'_i = \mu_i \tau^{}_i$.
For $\gamma \in \centralizer_{\Gamma} (\Lambda_v)$ let $(h^{}_1 , h^{}_2)$ and $(h'_1 , h'_2)$ denote the coordinates of the translation vectors $\gamma_{\EE}$ and $\gamma'_{\EE}$ with respect to the Basis $\{ \xi^{}_1 , \xi^{}_2 \}$.
For $i,j$ with $i \neq j$ we then have $\tau_i (\gamma) = < h_1 \xi_1 + h_2 \xi_2 , w_i > = h_j \tau_i (\xi_j)$ and $\tau'_i (\gamma) = < h'_1 \xi^{}_1 + h'_2 \xi^{}_2 , w'_i > = h'_j \tau'_i (\xi^{}_j)$.
The equality $h'_j \tau'_i (\xi^{}_j) = \tau'_i (\gamma) = \mu^{}_i \tau^{}_i (\gamma) = \mu^{}_i h^{}_j \tau^{}_i (\xi^{}_j) = h^{}_j \tau'_i (\xi^{}_j)$ then yields $h'_j = h^{}_j$ since $\tau'_i (\xi^{}_j) = \mu^{}_i \tau^{}_i (\xi^{}_j) \neq 0$.
Thus, after possibly changing the basis%
\footnote{Equality of the geometric data $\Tau$ and $\Tau'$ depends on neither the particular choice of the generators of $\Lambda_v$ nor the basis of $\EE_v$.}
of $\EE_v$, we have $\Tau^{}_v (\gamma) = (h^{}_1,h^{}_2) = (h'^{}_1,h'^{}_2) = \Tau'_v (\gamma)$.

The geometric data $\smash{\Tau^{}_{(\,.\,)}}$ and $\smash{\Tau'_{(\,.\,)}}$ therefore coincide at any vertex adjacent to consecutive edges $e,f$ with $\smash{\check{\T} (e) \neq \check{\T} (f)}$.

\providecommand{\bysame}{\leavevmode\hbox to3em{\hrulefill}\thinspace}
\providecommand{\MR}{\relax\ifhmode\unskip\space\fi MR }
\providecommand{\MRhref}[2]{%
  \href{http://www.ams.org/mathscinet-getitem?mr=#1}{#2}
}
\providecommand{\href}[2]{#2}

\vspace{4em}

\noindent
\textsc{Institute for Algebra and Geometry\\Karlruhe Institute of Technology, 76133 Karlsruhe (Germany)}

\medskip

\noindent
\texttt{\href{mailto:sebastian.grensing@kit.edu}{sebastian.grensing@kit.edu}}

\end{document}